\documentclass{article}

\usepackage[preprint]{neurips_2024}
\usepackage{xcolor}
\definecolor{myred}{RGB}{255, 0, 0}
\usepackage{amsbsy}
\usepackage{marvosym}
\usepackage{tikz-cd}
\usepackage{booktabs,array}

\usepackage[utf8]{inputenc} 
\usepackage[T1]{fontenc}    
\usepackage{hyperref}       
\usepackage{url}            
\usepackage{booktabs}       
\usepackage{amsfonts}       
\usepackage{nicefrac}       
\usepackage{microtype}      
\usepackage{xcolor}         
\usepackage{graphicx}
\usepackage{amsmath}
\usepackage{subfigure}

\newcommand{\R}{\mathbb{R}}

\title{PDEformer-1: A Foundation Model for One-Dimensional Partial Differential Equations}

\author{
  Zhanhong Ye\\
  Peking University\\
  Beijing, China \\
  \texttt{yezhanhong@pku.edu.cn} \\
 \And
   Xiang Huang \\
  Huawei Technologies Co. Ltd\\
  Hangzhou, China\\
  \texttt{huangxiang42@huawei.com} \\
  \And
  Leheng Chen \\
  Peking University \\
  Bejing, China\\
  \texttt{chenlh@pku.edu.cn} \\
    \And
    Zining Liu \\
    Peking University \\
    Beijing, China \\
    \texttt{liuzining@stu.pku.edu.cn} \\
    \And
    Bingyang Wu\\
    Peking University\\
    Beijing, China\\
    \texttt{wby2003@stu.pku.edu.cn}
    \And
    Hongsheng Liu\\
    Huawei Technologies Co. Ltd\\
    Hangzhou, China\\
    \texttt{liuhongsheng4@huawei.com}
    \And
    Zidong Wang\\
    Huawei Technologies Co. Ltd\\
    Hangzhou, China\\
    \texttt{wang1@huawei.com}
    \And
    Bin Dong\textsuperscript{\Letter} \\
    Peking University \\
    New Cornerstone Science Laboratory \\
    Beijing, China\\
    \texttt{dongbin@math.pku.edu.cn} \\
}

\begin{document}

\maketitle

\begin{abstract}
This paper introduces PDEformer-1, a versatile neural solver capable of simultaneously addressing various partial differential equations (PDEs).
With the PDE represented as a computational graph, we facilitate the seamless integration of symbolic and numeric information inherent in a PDE. 
A graph Transformer and an implicit neural representation (INR) are employed subsequently to generate mesh-free predicted solutions. 
We generated a dataset with up to three million samples involving diverse one-dimensional PDEs to pretrain our model.
Compared with baseline models trained specifically on benchmark datasets, our pretrained model achieves comparable accuracy via zero-shot inference, and the advantage expands after finetuning.
For PDEs new or unseen in the pretraining stage, our model can adapt quickly by finetuning on a relatively small set of examples from the target equation.
Additionally, PDEformer-1 demonstrates promising results in the inverse problem of PDE scalar coefficient recovery and coefficient field recovery.
\end{abstract}

\section{Introduction}
Partial differential equations (PDEs) are fundamental in describing the physical world.
They are widely used in various fields, such as physics, engineering, and finance. 
Solving PDEs is a challenging task due to the complex nature of the equation and the complex structure of the solution space. 
Traditional numerical methods, such as finite difference, finite element, and finite volume methods, have been developed to solve PDEs. 
However, these methods are computationally expensive. 
Moreover, they are unsuitable for real-time PDE solving and are not flexible enough to handle different PDEs without modification, which limits their application in practice. 
Recently, deep learning has shown great potential in solving PDEs. 
Neural networks can approximate the solution of PDEs with satisfactory accuracy and efficiency. 
Inspired by the great success of general-purpose models in computer vision (CV) and natural language processing (NLP), we construct a foundation of numerical PDE solving, named PDEformer-1, so that it can solve a wide range of PDEs.
Idealy, for an arbitrary PDE, our model can predict the solution with satisfactory accuracy and efficiency either with zero-shot inference or with inference after finetuning on a small number of samples. 

Compared with specialized models trained for specific PDEs, such as DeepONet\citep{lu_deeponet_2021} and FNO \citep{li_fourier_2021}, this foundation model can bring significant advantages.
Once trained, the model can solve various PDEs without retraining, saving a lot of time and resources, which differs from specialized models, as they need to be trained from-scratch for every new equation.
Furthermore, different forms of PDEs may share some common features.
For example, the same physical laws (diffusion, convection, reaction, etc.) can appear in various systems, and the influence of initial and boundary conditions on PDE solutions could be similar.
A general-purpose model can potentially capture and make use of these commonalities across different PDEs, 
and transfer the learned knowledge to a new PDE with scarce data to improve the prediction accuracy.
Indeed, constructing such a general-purpose model is not an easy task.
The model should be able to capture the essential features of PDEs and learn the mapping from the PDE to the solution, which may require complex network architecture design, massive data, advanced training strategies, and enormous computational resources. 
However, we hold the belief that such a general-purpose model can be well constructed and will significantly benefit the field of PDE solving.

In this paper, we present PDEformer-1, a versatile neural network designed to tackle a wide range of PDEs concurrently.
Our approach implements a unified model for one-dimensional PDEs, representing the PDE as a computational graph to seamlessly integrate the symbolic and numeric data. 
By incorporating flexible numeric encoders, a graph-level Transformer-based encoder, and an INR decoder, PDEformer-1 accepts input fields (initial conditions, source terms, etc.) given on arbitrary spatial points, and generates mesh-free predicted solutions.
After pretraining on a dataset containing a variety of PDEs, our model achieves high zero-shot accuracy on benchmark tests, outperforming well-trained specialized models. 
Furthermore, it can adapt quickly to new, unseen PDEs through finetuning on a smaller set of examples. 
PDEformer-1 also shows great potential in solving inverse problems such as PDE scalar coefficient recovery and coefficient field recovery.

A preliminary version of this work appears as~\citet{pdeformer}.
In this work, the network architecture is adapted so that the input coefficient fields can be defined on non-uniform grids.
The dataset we use covers a more diverse class of PDEs.
More experiments are conducted to evaluate the performance of our model.
Our code, pretrained model weights as well as the training dataset are available at \url{https://gitee.com/mindspore/mindscience/blob/master/MindFlow/applications/pdeformer1d}.

This paper is organized as follows.
Section~\ref{sec:related} briefly summarizes previous works using deep learning to solve PDEs.
Section~\ref{sec:methodology} introduces the architecture of our PDEformer-1 model.
Section~\ref{sec:forward} presents experimental results on forward problems, including zero-shot inference and inference after finetuning. Section~\ref{sec:inverse} presents experimental results on inverse problems, including scalar coefficient recovery, system identification, and coefficient field recovery.
Some further discussions are given in Section~\ref{sec:discussion}.

\section{Related Works}
\label{sec:related}
Generally, a PDE to be solved would involve two parts of information: one is the symbolic part, which specifies the mathematical form of the PDE, and the other is the numeric part, including coefficients, initial values, and boundary values of the PDE.
We categorize the existing deep learning-based PDE-solving methods into the following three classes based on their generality:

\paragraph{Independent Neural PDE Solvers for Specific PDE.}    In this category, a neural network is trained to solve a specific PDE whose symbolic representation and numeric information are fixed. Most models of this kind take a neural network as a surrogate of the solution, and the network is trained based on the information directly provided by the PDE, such as the initial value, boundary value, and the PDE itself. Unlike traditional numerical methods, these models do not require the discretization of the domain and can directly predict the solution at any point in the domain. 
For example, Physics-Informed Neural Networks (PINN) \citep{raissi_physics-informed_2019} and Deep Galerkin Method (DGM) \citep{sirignano_dgm_2018} take the equation and the boundary conditions as a penalty function to train the neural network. 
Deep Ritz Method (DRM) \citep{e_deep_2018} facilitates the variational form of the PDE to express the solution as a minimizer of a functional.
Weak Adversarial Networks (WAN) \citep{zang_weak_2020} use neural networks to represent the weak solution and trial function separately and train the network with a similar strategy as adversarial training.
The approaches mentioned above are all capable of solving specific PDEs without supervision, avoiding the need for domain discretization and training data generation.
However, these models are not general and need to be retrained for different PDEs and coefficients, which limits their ability to solve PDEs efficiently.

\paragraph{Parameterized Neural PDE Solvers.}    Compared to the previous category, models of this kind are designed to solve PDEs with fixed symbolic representation but variable numeric information.
These methods can be divided into two subcategories.
The first subcategory parameterizes the operator from coefficients to solutions using neural networks, typically referred to as neural operators.
Most neural operators require the use of traditional solvers to generate large datasets for a specific equation to conduct supervised learning, such as DeepONet \citep{lu_deeponet_2021} and FNO \citep{li_fourier_2021}. 
The notable advantage of these approaches is that, once training is completed, online solving of the same equation and similar parameters requires only a single neural network inference, with the time consumed being almost negligible. 
Moreover, the computational accuracy is sufficient to meet the demands of many real applications. 
The second subcategory can be seen as improvements of the single-PDE solving methods in the first category. 
These models are combined with meta-learning or transfer learning to solve parametrized PDEs, such as Meta-Auto-Decoder (MAD)~\citep{ye2023meta} and NRPINN~\citep{liu2022novel}, which exhibit higher accuracy when facing PDE coefficients that are not included in the training dataset (out-of-distribution, OoD). 
However, both subcategories require training a specialized model for each PDE, which brings great computational costs. 
It is also challenging to leverage the commonalities between equations of different forms to improve the performance of the neural solvers.

\paragraph{General-Purpose PDE Solvers}
Besides PDEformer-1, there already exist several attempts to construct a general-purpose PDE solver that can solve PDEs of different forms.
One major difficulty in building such a general-purpose model lies in how to integrate the symbolic information of PDEs efficiently.
Previous explorations are done mainly through two technical paths.
The first path is importing a language model into the model, taking the textual representation of the PDE as input, including Physics-Informed Token Transformer (PITT) \citep{lorsung_physics_2024}, PRedicting Operators and Symbolic Expressions (PROSE) \citep{liu_prose_2023}, and ICON-LM \citep{ICON-LM}.
These works separate symbolic information from numeric information and use different networks to process them, ignoring the complex relationship between these two parts.
The second path is to implicitly merge the symbolic information into the numeric information rather than explicitly taking it as part of the network input.
For example, inspired by in-context learning in the field of NLP, In-Context Operator Network (ICON) \citep{ICON} and its subsequent work \citep{ICON2} prepare in advance multiple sets of data pairs containing PDE parameters and solutions for the equations to be solved, in which different forms of equations will result in different data pairs.
These are taken as input into the network along with new parameters to be solved to generate predicted solutions.
Specialized in time-dependent PDEs, Multiple Physics Pretraining (MPP) \citep{MPP} and OmiArch \citep{OmniArch} input the solution trajectory history containing multiple time steps into the network and predict the solution for the next time step.
In this way, the information on the form of the equation can be reflected to a certain extent by these historical trajectories of solutions.
This type of implicit input method may not be sufficient to handle more diverse and complex types of PDEs, and often needs to work in conjunction with another solver to prepare additional solution data as network input, which limits their applicability.

\section{Method}
\label{sec:methodology}
\subsection{Overview of Our Model}
\begin{figure}[htpb]
	\centering
	\includegraphics[width=\linewidth]{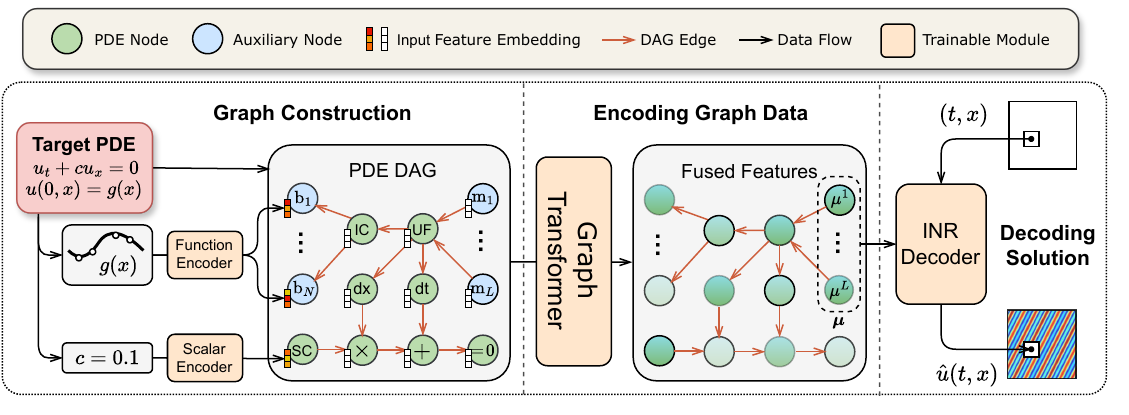}
	\caption{Overall architecture of PDEformer-1, taking the Advection equation $u_t+cu_x=0$, $u(0,x)=g(x)$ with periodic boundary condition as the example.
	Since only one unknown field variable $u_0=u$ is involved in this equation, we omit the subscript $j=0$ for $u_j$ and $\pmb{\mu}_j$ in the figure.}%
	\label{fig:PDEformerV2Arch}
\end{figure}

We consider one-dimensional time-dependent PDEs on $(t,x)\in [0,1] \times [-1,1]$ of the general form:
\[\mathcal{F}(u_1,u_2,\dots,c_1,c_2,\dots,s_1(x),s_2(x),\dots)=0,\]
where $c_1,c_2,\dots$ are real-valued coefficients, $s_1,s_2,\dots$ are scalar functions that may serve as initial conditions or coefficient fields in the equation, and $u_1,u_2,\dots$ are unknown field variables to be solved in the equation.
Here, we assume that operator $\mathcal{F}$ has a symbolic expression, which may involve differential and algebraic operations.
The goal is to construct a surrogate of the solution mapping $(\mathcal{F},c_1,c_2,\dots,s_1,s_2,\dots)\mapsto (u_1,u_2,\dots)$,
which takes both the PDE's symbolic expression $\mathcal{F}$ and the numeric information (including $c_1,c_2,\dots,s_1,s_2,\dots$) as model inputs and predicts the solution of the PDE ($u_1,u_2,\dots$).
We illustrate the overall network architecture in Figure~\ref{fig:PDEformerV2Arch}.
As shown in the figure, PDEformer-1 consists of three parts to solve the PDE: graph construction, encoding graph data, and decoding the solution.
In the following text, we shall elaborate on the details within these three parts.

\subsection{Graph Construction}
We first represent $\mathcal{F}$, i.e., the symbolic information specifying the PDE form, as a computational graph.
In such a computational graph, a node may stand for an unknown field variable (denoted as \verb|UF|), a scalar coefficient (\verb|SC|), an initial condition (\verb|IC|), a boundary value ($\mathtt{|x_L}$ or $\mathtt{|x_R}$), as well as a differential or algebraic operation.
A directed edge can be used to specify the operands involved in an operation.
These nodes and edges constitute a directed acyclic graph (DAG) with heterogeneous nodes and homogeneous edges.
Further details and examples of such a computational graph can be found in Appendix~\ref{sec:PDEDAGdetails}.

Then, in order to include the numeric information, we endow each graph node with an input feature embedding vector $\xi_i\in\R^{d_e}$.
For a scalar coefficient $c$, we input this numeric value into a scalar encoder and use the $d_e$-dimensional output as the input features for the corresponding \verb|SC| node.
The case for the boundary value nodes $\mathtt{|x_L}$ and $\mathtt{|x_R}$ is analogous.
In terms of a function $s(x)$ that typically contains richer information, we shall introduce $N$ new `branch' nodes, whose types are denoted as $\mathtt{b}_1,\mathtt{b}_2,\dots,\mathtt{b}_N$, respectively, and connect them to the \verb|CF| or \verb|IC| node corresponding to $s(x)$.
The function $s(x)$ is represented numerically as a series of scattered points $\{(x_j,s(x_j))\mid j\in \mathcal{J}\}$.
We feed these scattered points into a function encoder and split its $d_eN$-dimensional output vector into $N$ parts to form all the $d_e$-dimensional input features of those branch nodes.
The input features of all the remaining nodes are set to zero vectors. \footnote{To simplify implementation, these are actually the output of the scalar encoder when it receives a zero input.}
More implementation details can be found in Appendix~\ref{sec:fun_enc}.

Moreover, for each unknown field variable to be solved, we introduce $L$ additional nodes with types $\mathtt{m}_1,\mathtt{m}_2,\dots,\mathtt{m}_L$, respectively, and connect them to the corresponding \verb|UF| node.
These nodes will be used to decode the predicted solution, as will be explained below.

\subsection{Enocding Graph Data}
The graph data obtained in the previous step contains both symbolic and numeric information inherent in the PDE.
We shall integrate the information from the graph data and generate a latent code that represents the solution for each field variable $u_j$ to be solved. 
Each latent code takes the form 
\[\pmb{\mu}_j = [{\mu}^1_j, \dots, {\mu}^L_j]^{\mathrm{T}} \in\R^{L \times d_e}.\]

The fusion process leverages a graph Transformer, 
a powerful type of graph neural network based on the Transformer architecture, which is skilled at capturing and expressing complex graph structural information.
In practice, we adopt the Graphormer architecture~\citep{ying_transformers_2021} and make some adjustments to adapt it specifically for encoding PDEs. 
The details are presented in Appendix~\ref{sec:detail_DAG}.
For each variable index $j$ and layer index $\ell$, the \verb`UF` node representing $u_j$ is connected with another node with type $\mathtt{m}_\ell$, and we let $\mu_j^\ell\in\R^{d_e}$ be the embedding vector assigned to this $\mathtt{m}_\ell$ node in the output layer of the graph Transformer.

\subsection{Decoding the PDE solution}
We employ an INR that takes the coordinate $(t,x)$ as input and produces the mesh-free prediction $\hat{u}_j(t,x)$ according to $\pmb{\mu}_j$ for each field variable.
Various INR architectures with such an external condition have been adopted in neural operators~\citep{DINo}, data compression~\citep{coinpp}, and generative models~\citep{PolyINR}. 
In PDEformer-1, we utilize an adapted version of Poly-INR~\citep{singh2023polynomial} with $L$ hidden layers due to its efficiency,
and the modulations of the $\ell$-th hidden layer is generated according to $\mu_j^\ell$.
Details of the specific INR architecture are presented in Appendix~\ref{sec:inr}

\section{Results on Forward Problems}
\label{sec:forward}
\subsection{Pretraining Dataset}
We generate the dataset to pretrain our PDEformer-1 using diffusion-convection-reaction equations (with only one unknown field variable) of the following general form:
\begin{equation}\label{eq:pretrain}\begin{split}
  u_t+f_0(u)+s(x)+(f_1(u)-\kappa(x)u_x)_x&=0 , \quad (t,x) \in [0,1] \times [-1,1], \\
  u(0,x) &= g(x), \quad x \in [-1,1],
  \end{split}\end{equation}
where $f_i(u)=c_{i1}u+c_{i2}u^2+c_{i3}u^3$ for $i=0,1$.
Each coefficient $c_{ik}$ is set to zero with probability $0.5$, and drawn randomly from $\mathcal{U}([-3,3])$ otherwise. 
Note that terms with zero coefficients will be excluded in the computational graph of the PDE.
The coefficient fields $s(x),\kappa(x)$ are both selected from the following cases with certain probability: a random field, a random scalar (with no spatial dependency), and zero.
The value of $\kappa(x)$ is confined to the range $[10^{-3}, 1]$ except for the zero-valued case.
The initial condition $g(x)$ and the non-constant source term $s(x)$ are generated in the same way as the PDEBench dataset~\citep{PDEBench}. 
For the non-periodic case, the boundary condition at the left endpoint $x_L=-1$ takes the general form
\[(\alpha_Lu+\beta_Lu_x)|_{x=x_L}=\gamma_L.\]
Here, the type of the boundary condition is randomly selected from Dirichlet ($\alpha_L=1,\beta_L=0$), Neumann ($\alpha_L=0,\beta_L=1$) and Robin ($\alpha_L=\cos\theta_L,\beta_L=\sin\theta_L,\theta_L\sim\mathcal{U}([0,\pi])$) with equal probability.
The specific boundary value $\gamma_L$ is randomly selected from zero, a random scalar drawn from $\mathcal{U}([-3,3])$, and the initialized value $\alpha_Lg+\beta_Lg_x$ with equal probability.
The boundary condition at the right endpoint $x_R=1$ is determined analogously and independently.

\begin{figure}[htbp]
\centering
\includegraphics[width=0.6\linewidth]{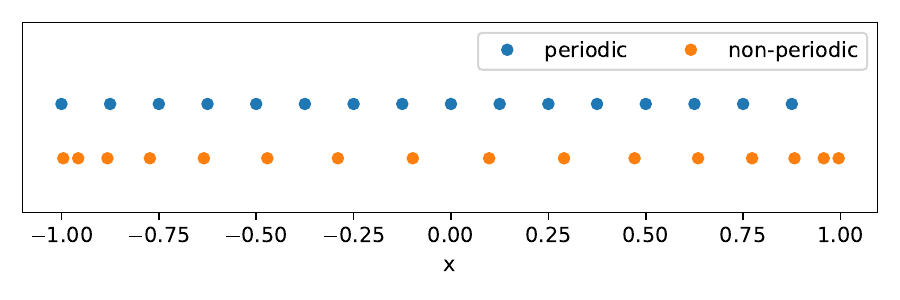}
\caption{An illustrative comparison of the grids for periodic and non-periodic PDEs in our dataset, in which a resolution of $16$ points is used for simplicity.
Thanks to the mesh-free nature of PDEformer-1, we can utilize solution samples recorded on different grid points directly during training, and no interpolation is required.}%
\label{fig:dedalus_grid}
\end{figure}

We use the open-source software package, Dedalus V3~\citep{Dedalus}, to generate numerical solutions and compose the training dataset. 
Based on spectral methods, the solver employs a uniform grid in the periodic case, while in the non-periodic case, the grid points cluster quadratically towards the endpoints of the spatial interval, as can be seen in Figure~\ref{fig:dedalus_grid}.
In both cases, we discretize the spatial solution domain into $256$ grid points.
The solver proceeds at a time-step of $\delta t_\text{solver}=4\times 10^{-4}$, and the solution snapshots are recorded with time-step $\delta t_\text{data}=0.01$, yielding a total of $101$ temporal values for each data sample.
If Dedalus fails to provide a solution for the PDE, or if the generated solution does not satisfy $\|u\|_{L^\infty}\le 10$, we shall discard the corresponding data sample, and not include it in the final dataset.

\subsection{Scaling Law}
\begin{figure}[ht]
    \centering
    \includegraphics[width=\linewidth]{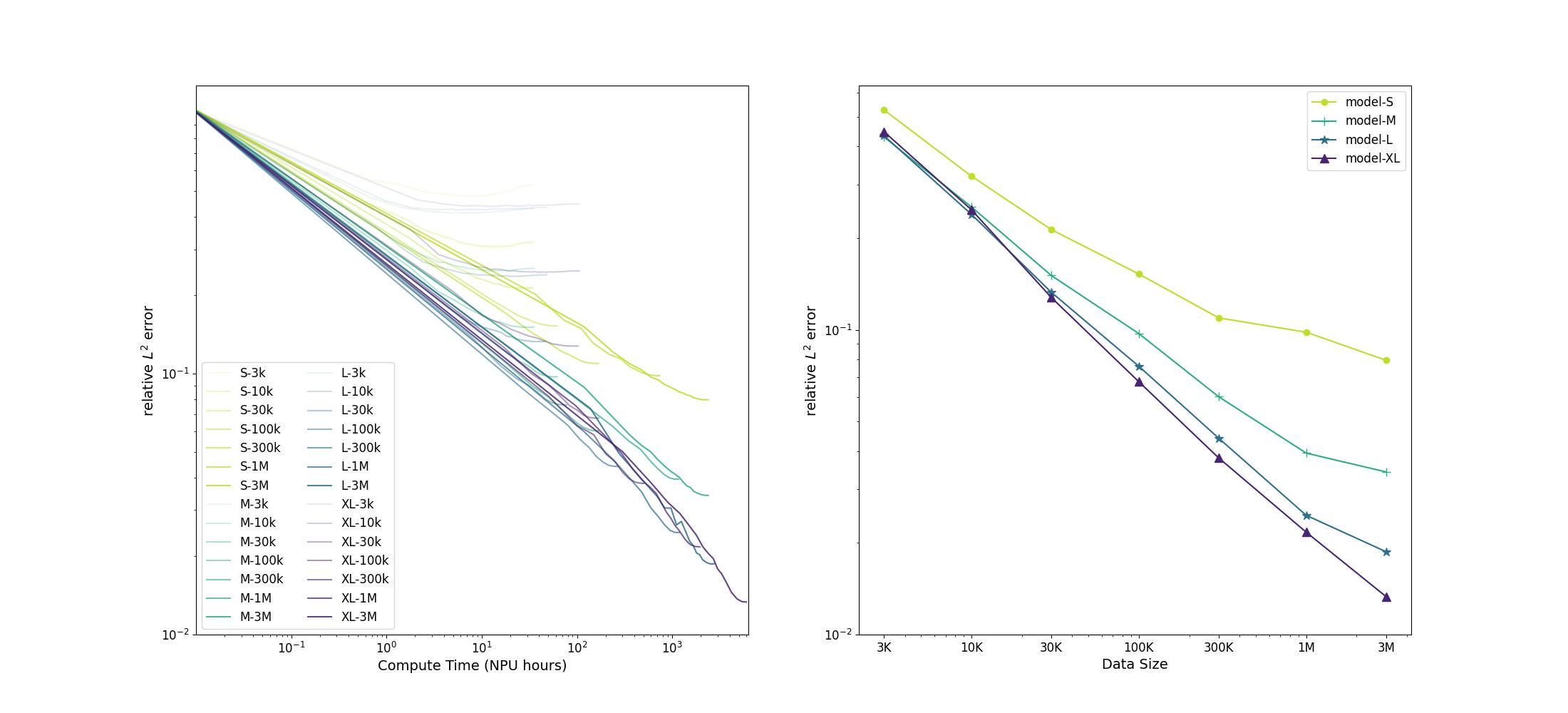}
    \caption{The figure on the left shows the change of test accuracy according to compute time. The figure on the right shows the change of test accuracy according to training dataset size. Different lines refers to different scales of model or training dataset. For example, M-30k represents training the M-size model on a dataset with 30k data samples.}%
    \label{fig:4}
\end{figure}

We conduct experiments on PDEformer-1 with various scales of model (S, M, L, and XL)\footnote{PDEformer-1 of scale S, M, L, and XL differs mainly in the depth and width of the graph Transformer and INR. The number of parameters of models of scale S, M, L, and XL are 3.60M, 6.92M, 22.4M, and 58.2M, respectively. We elaborate further details of the model in Appendix~\ref{sec:setting}.} and dataset constructed the same as the pretraining dataset, and the results of training time and test accuracy are shown in Figure~\ref{fig:4}.  
From the figure, it can be observed that with a fixed number of epochs, the test error of the model exhibits a power-law relationship with the amount of data. This relationship can be mathematically expressed as $\text{Test-Error} \propto \text{Data-Amount}^{-\alpha}$.
Additionally, in the case where the training data and model parameters are sufficiently large, the test error exhibits a power-law relationship with the training time, which can be written as $\text{Test-Error} \propto \text{Training-Time}^{-\beta}$.

\subsection{Performance on In-Distribution Test Sets}
\begin{figure}[ht]
\centering
\includegraphics[width=\linewidth]{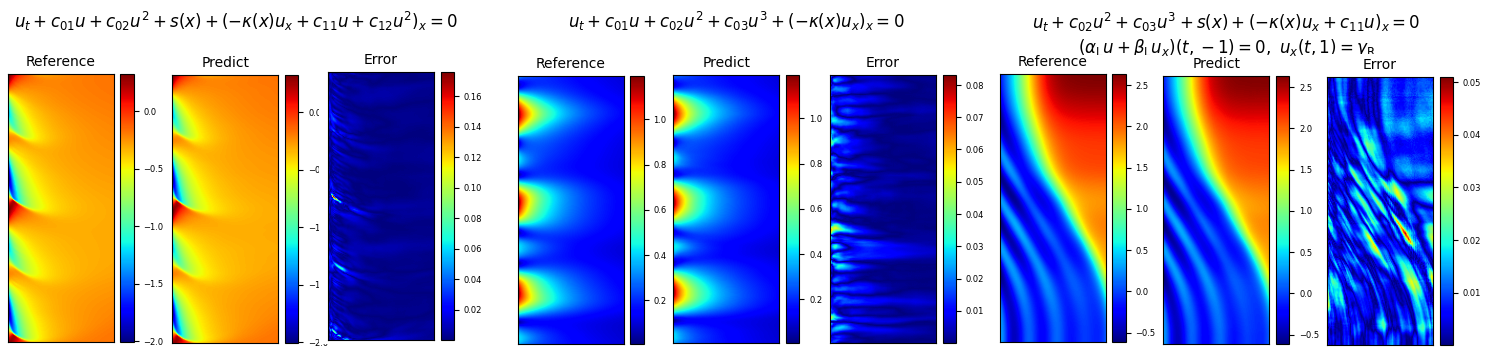}
\caption{Inference performance of PDEformer-1 on part of the in-distribution test set. The horizontal axis represents the temporal axis $t$, and the vertical axis represents the spatial axis $x$.}%
\label{fig:47}
\end{figure}

After pretraining PDEformer-1 (size XL)\footnote{All the following experimental results are conducted with our XL-sized model, unless specified otherwise.} for 125 hours on 48 NPUs with 3 million data samples, the model achieves a normalized root-mean-squared-error (nRMSE) of $0.0084$ on the pretraining dataset.
To test its robustness, we generate a test dataset in the same manner as the pretraining dataset (which we refer to as in-distribution), including 10k periodic samples and 10k non-periodic samples.
PDEformer-1 achieves an nRMSE of $0.0133$ on this test set.
Figure~\ref{fig:47} shows the prediction results of PDEformer-1 on the test set.
From the results, it can be seen that for equations distributed same as the pretraining dataset, our PDEformer-1 model achieves satisfactory accuracy in solving equations.

\paragraph{PDEBench Dataset}
PDEBench~\citep{PDEBench} includes Burgers', Advection, and 1D Reaction-Diffusion PDEs.
We shall examine some scenarios where all the PDE parameters lie within the distribution of the pretraining data, including Burgers' equation with $\nu=0.1$ and $\nu=0.01$, as well as the Advection equation with $\beta=0.1$.
We select U-Net~\citep{U-Net}, FNO~\citep{FNO}, and DeepONet~\citep{DeepONet} as benchmark models, performing direct inference with PDEformer-1 (Zero-Shot) and post-finetuning inference results with PDEformer-1 (Finetune), using the same datasets as the benchmark models.  Notably, all previous methods as well as PDEformer-1 (From-Scratch) are trained from-scratch and tested separately on different datasets.
The results are shown in Table~\ref{tab:error-id}.
\begin{table}[ht]
	\centering
	\begin{tabular}{cccc}
		\multicolumn{1}{c}{Model}  &\multicolumn{2}{c}{Burgers'}   & Advection    \\
									&\multicolumn{1}{c}{$\nu=0.1$} & \multicolumn{1}{c}{$\nu=0.01$} & \multicolumn{1}{c}{$\beta=0.1$} \\
		\hline \\
		U-Net & 0.1548            & 0.2222                     & 0.0897            \\
		DeepONet & 0.0687             & 0.1758                         & 0.0194  \\
		FNO & 0.0171             & 0.0492 & 0.0125              \\
        PDEformer-1 (From-Scratch) & 0.0383 & 0.0843 & 0.0413   \\
		PDEformer-1 (Zero-Shot)             & 0.0061            & 0.0180       & 0.0120    \\
		PDEformer-1 (Finetune)       & \textbf{0.0034}    & \textbf{0.0119}    & \textbf{0.0062}    \\
    \hline
	\end{tabular}
    \caption{Test relative $L_2$ error on PDEBench, in which the PDE coefficients lie within the range of the pretraining data.
    We format the best outcomes in bold.}
	\label{tab:error-id}
\end{table}
Remarkably, for Burgers' equation with $\nu=0.1$ and $0.01$, the zero-shot PDEformer-1 outperforms all the baseline models trained specifically on these datasets.
Such superior performance can be mostly attributed to the exposure to diverse PDEs during pretraining rather than the network architecture, as PDEformer-1 exhibits much more competitive performance after pretraining on a general dataset compared to training from-scratch.

\subsection{Inference Speed}
Table~\ref{tab:inference_time} showcases a comparison of the number of parameters, per-sample inference time and prediction accuracy for a range of models, including DeepONet, FNO, U-Net, and PDEformer.
We also include the results of two traditional numerical solvers.
The former is based on the first-order upwind finite-difference (FD) scheme, utilizing the \verb|solve_ivp| function provided by the SciPy Python package, and the latter being Dedalus, the spectral-method-based solver employed in generating our ground-truth solution data.
The evaluation was conducted using the 1D Advection equation ($\beta = 1.0$) on a $256 \times 256$ spatial-temporal grid as a test case, with neural network models tested on a single NPU and traditional solvers executed on a CPU.
The neural network models are adequately trained on the corresponding dataset, and the batch size is standardized to 10 during the test.
We average the total time consumption of each model across all samples to show the per-sample inference time.
As the FD solver exhibits lower accuracy, the spatial grid resolution is refined to $16 \times 256 = 4096$ in its solution process.
\begin{table}[ht]
\small
\centering
\begin{tabular}{lcccc|cc}
\textbf{Model} & \textbf{DeepONet} & \textbf{FNO} & \textbf{U-Net} & \textbf{PDEformer-1} & \textbf{FD} & \textbf{Dedalus} \\ \hline
\textbf{Num. Param.} & 1.65M & 0.92M & 13.39M & 58.24M & - & - \\
\textbf{Infer. Time (ms)} & 7.41 & 2.74 & 5.17 & 23.70 & 2072.3 & 410.8 \\
\textbf{Rel. $L_2$ Error} & 0.0194 & 0.0125 & 0.0897 & 0.0120 & 0.0674 & - \\
\end{tabular}
\caption{Comparison of the number of trainable parameters and per-sample inference time among models, where FD represents the traditional finite difference method, and Dedalus is the spectral method solver used to generate our pretraining dataset.}
\label{tab:inference_time}
\end{table}

\subsection{Performance of Finetuned Model on Out-of-Distribution Test Sets}
After pretraining, our PDEformer-1 is highly versatile in solving various equations.
To further evaluate its performance, we conduct experiments on zero-shot inference and post-finetuning inference on different datasets containing equations or coefficients not included in the pretraining dataset (i.e. OoD case).
Our evaluation involves benchmark PDEBench datasets as well as additional custom datasets generated by Dedalus.
\paragraph{PDEBench Dataset} 
In this scenario, we test our model on PDEs with coefficients (numeric information) lying outside the distribution range of the pretraining data, including Burgers' Equation with $\nu=0.001$, Advection Equation with $\beta=1$, and 1D Reaction-Diffusion Equation with $\nu=1, \rho=1$.
We finetune our model on the corresponding PDEBench dataset to adapt it to these OoD equations, and present the results in Table~\ref{tab:error-ood}.
In addition to the in-distribution experiments, the finetuned PDEformer-1 consistently excels in all out-of-distribution tests as well, further highlighting the robustness and versatility of our approach in solving a wide range of PDEs.
\begin{table}[ht]
	\centering
	\begin{tabular}{cccc}
		\multicolumn{1}{c}{Model}  & Burgers' & Advection & \quad Reaction-Diffusion \quad   \\
									&\multicolumn{1}{c}{$\nu=0.001$} & \multicolumn{1}{c}{$\beta=1.0$} & \multicolumn{1}{c}{$\nu=1,\rho=1$} \\
		\hline \\
		U-Net                    & 0.2404             & 0.2627          & 0.0144            \\
		DeepONet                 & 0.1988             & 0.0216          & 0.0019  \\
		FNO                      & 0.0727             & 0.0152          & 0.0021              \\
        PDEformer-1 (From-Scratch)  & 0.1078 &  0.0407 & 0.0027\\
		PDEformer-1 (Zero-Shot)         & 0.0388             & 0.9036          & 0.0377    \\
		PDEformer-1 (Finetune)      & \textbf{0.0265}    & \textbf{0.0099} & \textbf{0.0009}    \\
    \hline
	\end{tabular}
    \caption{Test relative $L_2$ error on PDEBench, in which the PDE coefficients lie outside the range of the pretraining data.
    We format the best outcomes in bold.}
	\label{tab:error-ood}
\end{table}

\begin{figure}[ht]
    \centering
    \includegraphics[width=\linewidth]{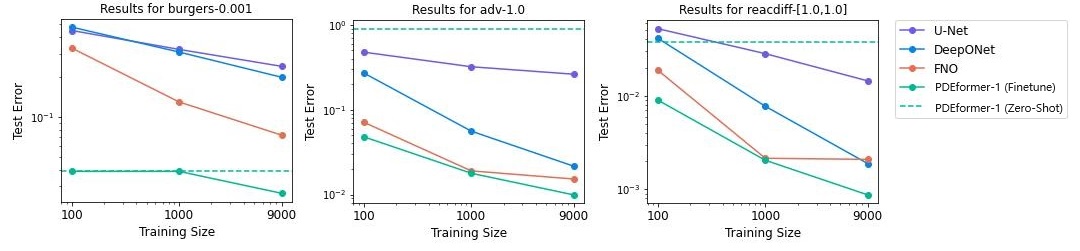}
    \caption{Variation of test error with number of finetuning samples. `PDEformer-1 (Zero-Shot)' represents our model's direct inference capability without any finetuning. This unique characteristic is visually depicted as a horizontal dashed line across the figures.}%
    \label{fig:5}
\end{figure}
We also embarked on a detailed investigation to assess the model's learning efficiency with limited data.
Specifically, we reduced the training dataset size from 9k to 1k and 100 samples.
As depicted in Figure \ref{fig:5}, the finetuned PDEformer-1 model notably excels, outperforming all the other methods in the test.
Moreover, for the case of Burgers' equation with $\nu=0.001$, the zero-shot PDEformer-1 establishes a commendably high benchmark, demonstrating robust performance without any finetuning.
It is particularly noteworthy that under OoD conditions, such as in the Advection $(\beta=1)$ and Reaction-Diffusion scenarios, the finetuned PDEformer-1 efficiently attains superior results.
This highlights the model's few-shot learning\footnote{In this work, we use the term \emph{few-shot learning} to describe the model's proficiency in adapting to and learning from new data that falls outside the distribution of the training set, using only a small number of examples for finetuning.} ability in adapting to unfamiliar scenarios. 

\begin{figure}[ht]
  \centering
  \includegraphics[width = 0.7\linewidth]{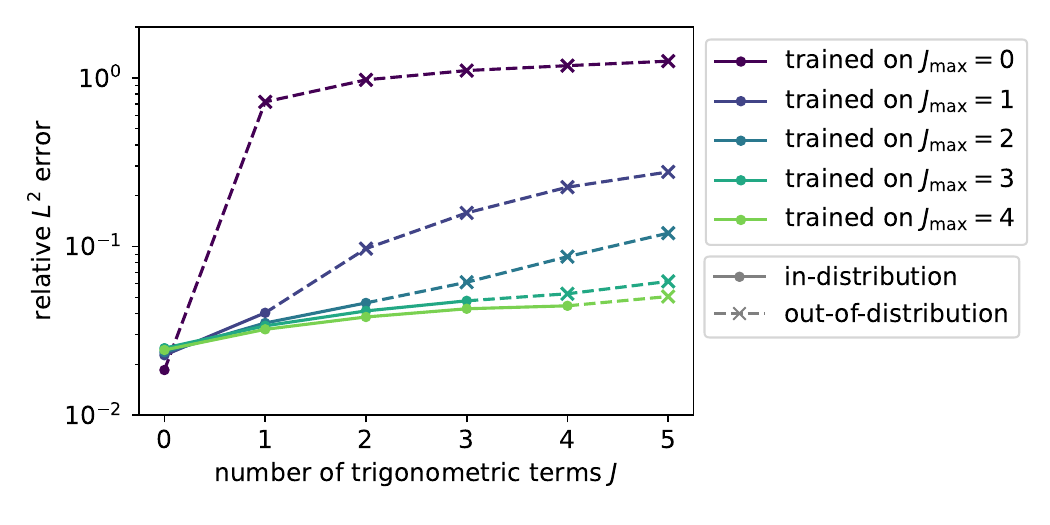}
  \caption{Prediction accuracy on PDEs with different number of trigonometric terms after finetuning on different datasets.}
  \label{fig:5(1)}
\end{figure}

\paragraph{Datasets with Trigonometric Functions}
In this scenario, the PDE forms (symbolic information) lie outside pretraining data distribution, involving trigonometric terms that is new to our model.
We consider PDEs with $J$ trigonometric functions of the general form
\[
    u_t + f_0(u) + s(x) + (f_1(u)-\kappa (x)u_x)_x = 0,\qquad (t,x)\in [0,1]\times [-1,1]  
,\]
where
\[f_i(u) = \sum_{k=1}^3 c_{i0k}u^k + \sum_{j=1}^{J_i}c_{ij0}h_{ij}(c_{ij1}u + c_{ij2}u^2),\]
each trigonometric function $h_{ij}\in\{\sin,\cos\}$ is randomly selected, and $J_0+J_1=J$ holds with $J_0 \in \{0,1,\dots,J\}$ randomly chosen as well.
Only periodic boundary conditions are employed for simplicity.
Note that when $J=0$, the PDE degenerates to the case we have used in the pretraining dataset~\eqref{eq:pretrain}.
For each $J\in\{1,2,3,4,5\}$, we generate a dataset of such equations containing 110k samples, among which 100k are used for training, leaving the rest 10k for testing.
For example, the case $J_{\max}=2$ means the pretrained PDEformer-1 model is finetuned on datasets corresponding to $J=0,1,2$ with a total of 300k data samples for $200$ epochs, and we then evaluate its performance on test datasets corresponding to $J=0,1,\dots,5$.
The number of finetuning epochs is adjusted for different values of $J_{\max}$ so as to ensure an equal number of iterations.
The results are shown in Figure~\ref{fig:5(1)}. 
According to the curve on the top of the plot, for the case $J_{\max}=0$, the model fails to solve PDEs with at least one trigonometric function, which is reasonable as no such knowledge about the trigonometric terms can be acquired during its training process.
As $J_{\max}$ increases, the prediction accuracy of the finetuned model increases accordingly.
It is noteworthy that, for the case $J_{\max}=4$, the finetuned model has never seen any data sample of PDEs with up to five trigonometric terms during its training process.
However, it still reaches a zero-shot prediction accuracy of $5.05\%$ on the $J=5$ test dataset.
This indicates that, after training on data samples with different numbers of trigonometric terms, our model has learned the semantic meaning of the summation operation applied on such terms, which enables its zero-shot generalization to more complicated PDEs that lie outside the distribution of the training dataset.

\paragraph{Wave Equation Dataset} 
\label{sec:wave}
Since the pretraining dataset does not contain second-order time derivatives, we finetune PDEformer-1 on the wave equation dataset to test its adaptability to new operators and generalization capabilities.
The wave equation takes the general form:
\begin{equation}\label{eq:wave}\begin{aligned}
	u_{tt}+\mu u_t+Lu+bu_x+f(u)+s_T(t)s_X(x)&=0,\\
	u(0,x)=g(x),\quad u_t(0,x)&=h(x),
\end{aligned}\end{equation}
where $f(u) = c_{1}u+c_{2}u^2+c_{3}u^3$, each coefficient $c_i$ is set to zero with probability 0.5, and drawn from the distribution $\mathcal{U}([-3,3])$ otherwise.
The second-order spatial term $Lu$ is randomly chosen from the following three forms with equal probability:
\[
L(x)=-c(x)^2u_{xx},\qquad
L(x)=-c(x)(c(x)u_x)_x,\qquad
L(x)=-(c(x)^2u_x)_x.
\]
The source term $s_T(t)s_X(x)$ is randomly chosen from the following five specific scenarios: zero, a random scalar, a random field (with only spatial dependency), a randomized varying scalar (with only temporal dependency), and a varying random field with spatiotemporal dependency.
Boundary conditions include periodic and non-periodic boundaries.
In the non-periodic case, the boundary conditions at both endpoints are chosen independently from Dirichlet, Neumann, Robin, and Mur ($(u_t-c(x)u_x)|_{x=x_L}=\gamma_L$ on the left, $(u_t+c(x)u_x)|_{x=x_R}=\gamma_R$ on the right) types with equal probability.

During the finetuning phase, we use 200k wave equation data samples each for periodic and non-periodic boundaries, along with the same number of data samples from the pretraining data.
The final nRMSE achieved is $5.1\%$ for periodic boundaries and $4.6\%$ for non-periodic boundaries, indicating that our model exhibits satisfactory performance on PDEs with new operators after finetuning.

\subsection{Efficiency of Finetuning}
\begin{figure}[htpb]
    \centering
    \includegraphics[width=\linewidth]{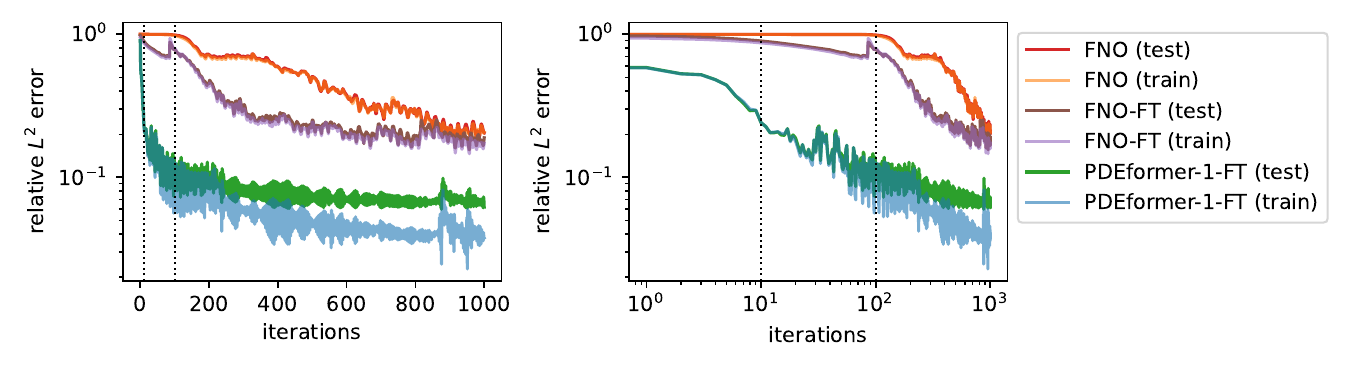}
    \caption{Comparison of the speed of finetuning PDEformer-1 with training FNO from-scratch and finetuning a pretrained FNO (FNO-FT) on the PDEBench Advection dataset with $\beta=1.0$.
    The right subfigure uses a logarithmic scale for the horizontal axis, whereas the left employs a linear scale.
    The vertical lines correspond to $10$ and $100$ iterations.
    }
    \label{fig:50}
\end{figure}

To compare the efficiency of finetuning PDEformer-1 with training traditional specialized models, we conduct a comparative analysis on the Advection equation ($\beta = 1$, OoD) dataset from PDEBench with a limited number of $100$ training samples.
As depicted in Figure~\ref{fig:50}, PDEformer-1 rapidly approaches convergence in about just 100 iterations.
Conversely, the FNO model, trained from-scratch, results in a higher test error even after one thousand iterations.
Indeed, it is possible for traditional neural operators to start from a better initialization.
However, designed for a specific type of PDE, they cannot be pretrained on 3M data samples containing diverse PDEs as PDEformer-1 does.
The valid option left for us is to pretrain them on one different PDE, and then transfer to the target setting, which could be much less efficient.
Post pretraining on 9k samples of the Advection equation with $\beta=0.1$ for 1k iterations, the \texttt{FNO-FT} model only exhibits a limited improvement over the corresponding from-scratch version, as can be seen in the figure.
This contrast highlights the pretrained PDEformer-1's swift and accurate adaptability, marking a significant advancement over existing specialized models.

\section{Results on Inverse Problem}
\label{sec:inverse}
To verify that our pretrained PDEformer-1 is applicable to downstream tasks, we conduct experiments on inverse problems as the example.
The tasks include recovering unknown scalar coefficients in equations, system identification, restoring source fields in diffusion-convection-reaction equations, and estimating velocity fields in wave equations. 
In these scenarios, we can obtain the information of PDE solutions under multiple initial conditions.
Instead of acquiring the full exact solution fields, the model only have access to the noisy observations at several spatial-temporal points of the solution, and the initial conditions available are also noisy.
We input the equation form, noisy initial conditions, and the current estimation (scalar coefficient or coefficient fields) into PDEformer-1, compare the predicted solutions with observations, calculate the relative $L_2$ error, and then solve the corresponding optimization problem to find the best estimations.

\subsection{Recovering Scalar Coefficients}
\begin{figure}[htbp]
  \centering
  \subfigure[][Results of recovered coefficients at noise level $3\%$\label{fig:15a}]{\includegraphics[width=0.35\linewidth]{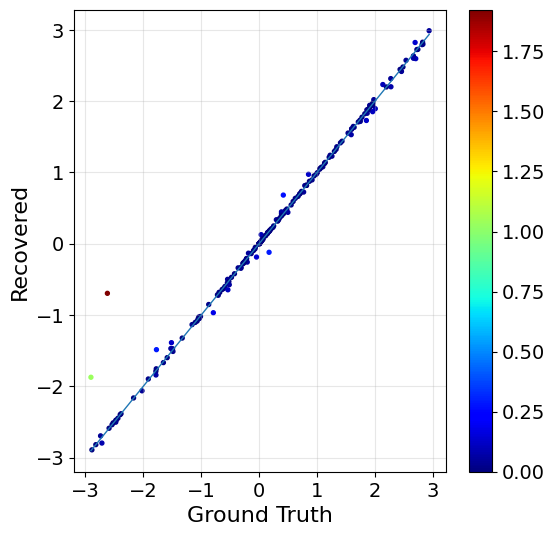}}
  \subfigure[][Noisy initial condition compared to ground truth\label{fig:15b}]{\includegraphics[width=0.55\linewidth]{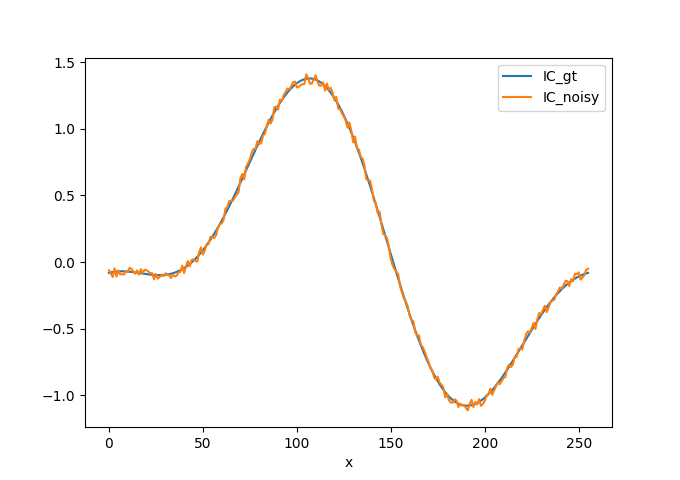}}
  \caption{Recovered results and example initial condition from the basic noise level}
  \label{fig:15}
\end{figure}

In the presence of known equation forms and noisy initial conditions, we recover the unknown coefficients in the equations based on noisy observations at certain spatiotemporal locations.
The forms of the equations used are consistent with those in the pretraining dataset, with periodic boundary conditions employed for simplicity. 
For example, given the equation $u_t + c_{11}u + c_{12}u^2 + s(x) + (-\kappa u_x)_x = 0$, we recover the scalar coefficients $c_{11},c_{12},\kappa$ from observations and known initial conditions. 
For each PDE, we assume that solutions under 8 different initial conditions can be obtained. 
The observable sample points are limited to 8 fixed spatial locations, and the accessible timesteps at each location are also randomly generated, with an average of 8 observation samples per spatial location.
Observational values at these points will have noise added with a relative magnitude of $10\%$, and the initial values that are inputs into PDEformer-1 will also have $3\%$ additive noise.
Due to the highly non-convex nature of the loss function for coefficients recovery, we employ the Particle Swarm Optimization (PSO)~\citep{kennedy_particle_1995} algorithm to solve the corresponding optimization problem. 
As can be seen from the results in Figure~\ref{fig:15}, even with only a limited number of noisy observations at a few locations, we can still recover most of the coefficient values in the equation successfully.

\begin{figure}
  \centering
  \subfigure[Nosie level $10\%$]{\includegraphics[width = \linewidth]{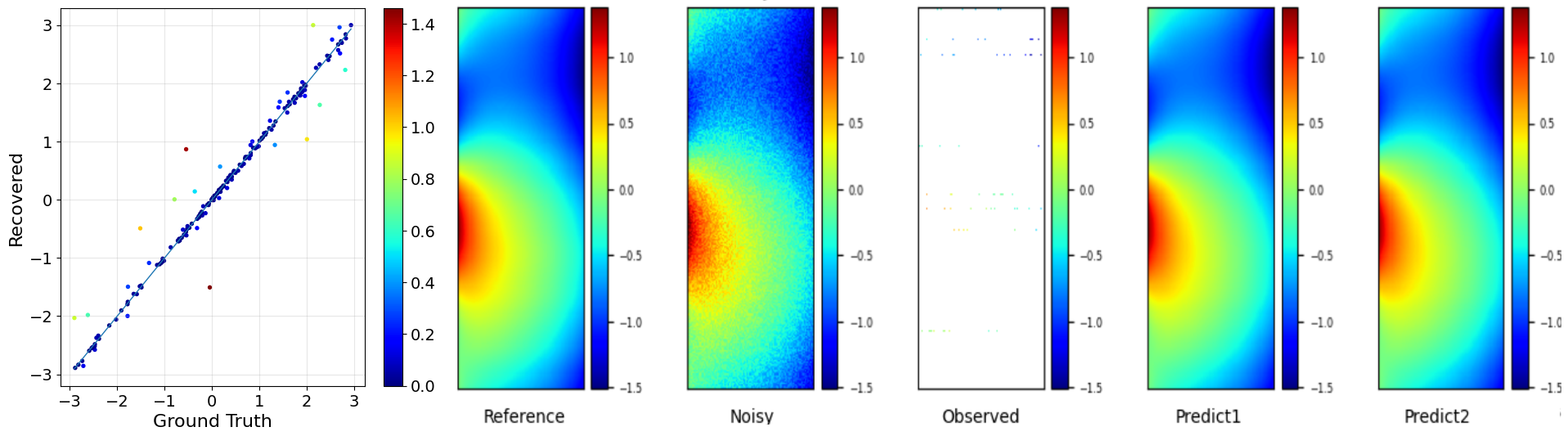}}
  \subfigure[Nosie level $30\%$]{\includegraphics[width = \linewidth]{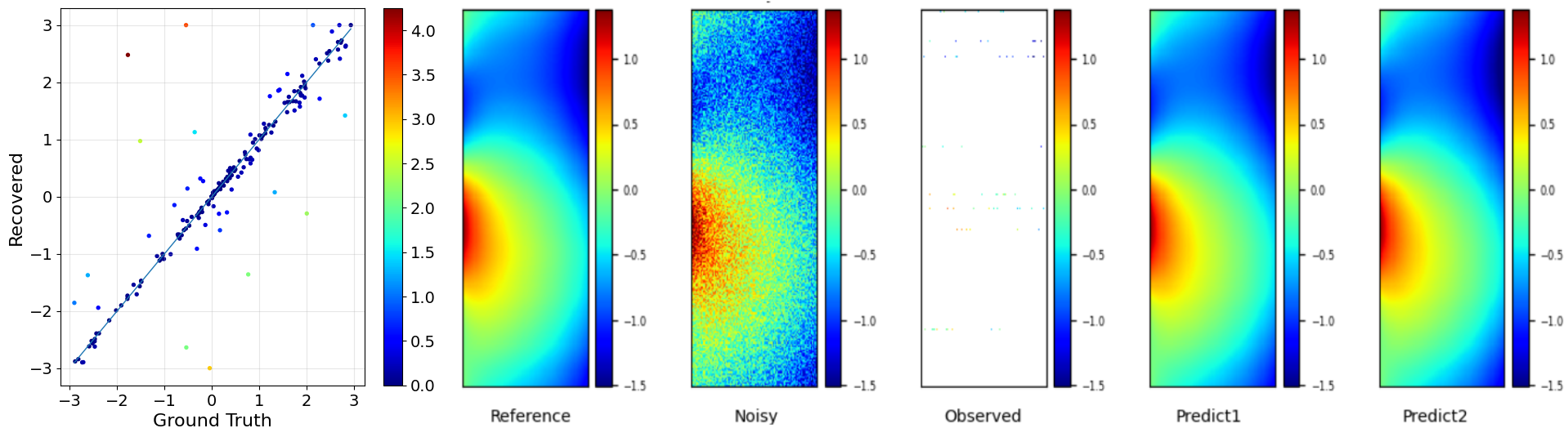}}
  \caption{Results of scalar coefficient recovery under different noise levels.
  The first column shows the results of our experiments using scattered points, in which a point lying on the diagonal line indicates perfect recovery of the corresponding coefficient.
  The remaining five columns, from left to right, contain the ground truth solution fields, solutions with additive noise, the sub-sampled noisy solutions that serve as the observation for the inverse problem, PDEformer-1 predictions using the ground truth coefficients, and predictions using the recovered coefficients, respectively.}
  \label{fig:scalar}
\end{figure}

We also conduct experiments with different noise levels, with the results shown in Figure~\ref{fig:scalar}.
Despite the increase in noise level, our PDEformer-1 performs well in recovering scalar coefficients, as the scattered points are still close to the diagonal line $y=x$.

\subsection{System Identification}
System identification is a methodological approach used in control engineering and signal processing to create mathematical models of dynamic systems based on measured data.
We conduct experiments on equation
\begin{equation*}
    u_t +c_{01}u + c_{02}u^2 + c_{03}u^3 + s(x) + \left( -\kappa u_x + c_{11} u + c_{12}u^2 + c_{13}u^3 \right)_x = 0
\end{equation*}
(same as~\eqref{eq:pretrain}), recovering coefficients\footnote{We treat the system identification problem as a special case of coefficient recovery. The only difference is that we do not assume prior knowledge about which coefficients are zero.}%
$c_{01},c_{02},c_{03},\kappa,c_{11},c_{12},c_{13}$ to determine the dynamic system behind the measured data. 
We assume that data can be obtained under 16 different initial conditions, with observable samples limited to 16 fixed spatial locations and an average of 16 observation samples per location.
We add noise with a relative amplitude of $1\%$ to the observed values at these points and $3\%$ noise to the initial values used as inputs for PDEformer-1. 
Part of the results are shown in table~\ref{table:sys_id}.
\begin{table}[ht]
\centering
\begin{tabular}{c|ccccccccc}
\toprule
 & $\kappa$ & $c_{01}$ & $c_{02}$ & $c_{03}$ & $c_{11}$ & $c_{12}$ & $c_{13}$ \\
\midrule
GT & 0.047 & 0.0 & 0.0 & 1.133 & 0.0 & -0.276 & 0.0 \\
Recover & 0.047 & -0.009 & 0.006 & 1.177 & -0.004 & -0.278 & 0.004 \\
\midrule
GT & 0.552 & 2.196 & 2.482 & 0.0 & -2.472 & 0.0 & 1.403 \\
Recover & 0.558 & 2.186 & 2.451 & 0.032 & -2.450 & 0.006 & 1.431 \\
\midrule
GT & 0.699 & 1.625 & 1.225 & -0.516 & 1.745 & -0.315 & 0.0 \\
Recover & 0.697 & 1.632 & 1.225 & -0.544 & 1.760 & -0.293 & -0.021 \\
\bottomrule
\end{tabular}
\caption{Results of system identification using PSO and PDEformer-1 as the surrogate of the forward model, where `GT' and `Recover' represent the value of ground truth and the recovered results by PSO and PDEformer-1. We can tell from the results that PDEformer-1 performs well in the problem of system identification and is robust to low-level noise.}
\label{table:sys_id}
\end{table}

\subsection{Restoring Source Fields} 
\begin{figure}
  \centering
  \includegraphics[width=\linewidth]{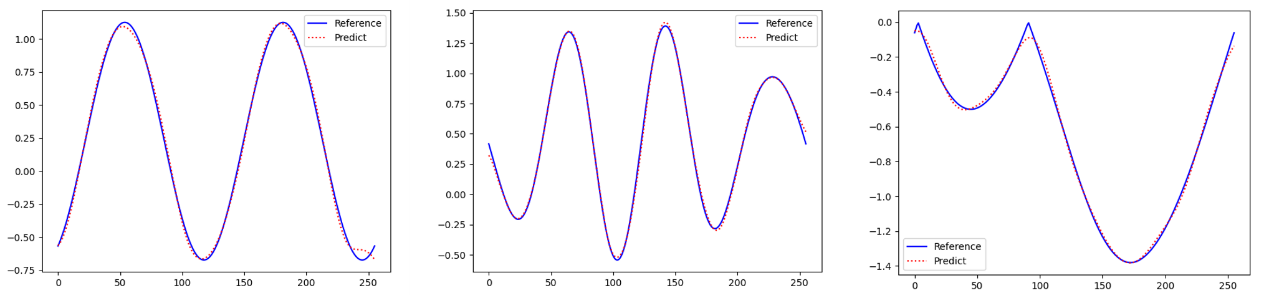}
  \caption{The results of the source term (coefficient field) recovery are depicted, where the red dashed lines represent the recovered source term, and the blue solid lines indicate its true value.}
  \label{fig:68}
\end{figure}

In this task, the form of the equations used still aligns with those in the pretraining dataset, and the objective is to recover the coefficient field $s(x)$ corresponding to the source term of the equations. For each PDE, we assume that solutions can be obtained under 100 different initial conditions, with observable sample points limited to 32 fixed spatial locations and an average of 32 observation samples available per location. Observational values at these points will have noise added with a relative amplitude of $1\%$, while the initial values input into PDEformer-1 will have $3\%$ noise added. The corresponding optimization problem is solved using the gradient descent method. From the experimental results shown in Figure~\ref{fig:68}, PDEformer-1 effectively accomplishes the task of recovering the source term in this inverse problem. 

\subsection{Recovering Velocity Fields in Wave Equations}
\begin{figure}[htpb]
	\centering
	\subfigure[$Lu = -(c^2(x)u_x)_x$]{\includegraphics[width=0.3\linewidth]{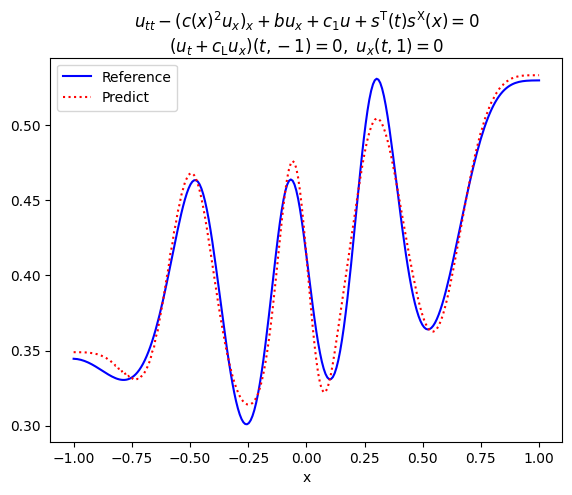}}
  \subfigure[$Lu = -c^2(x)u_{xx}$]{\includegraphics[width=0.3\linewidth]{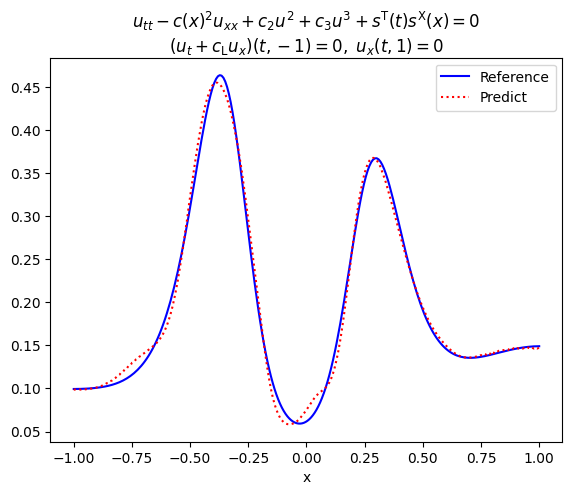}}
  \subfigure[$Lu = -c(x)(c(x)u_x)_x$]{\includegraphics[width=0.3\linewidth]{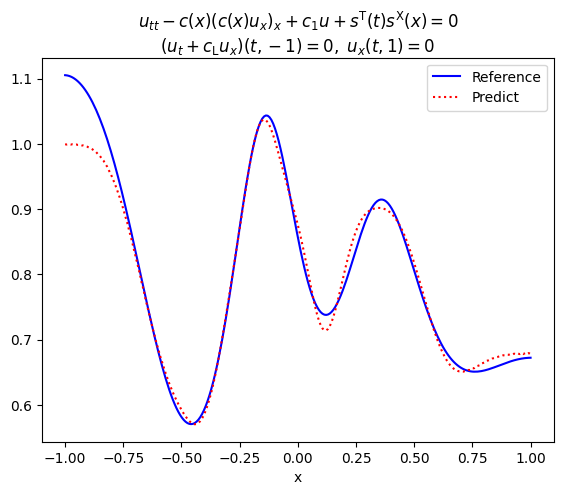}}
  \caption{Results of velocity field recovery, where the predicted coefficient field is plotted with red dashed lines, and the true value is plotted with blue solid lines.}
  \label{fig:70}
\end{figure}

After incremental training on wave equation data using PDEformer-1, as explained in Section~\ref{sec:wave}, we employ it to recover the unknown velocity field $c(x)$ in the wave equation.
To simulate the real-world full wave inversion scenario, we stick to source terms with spatial-temporal dependency in~\eqref{eq:wave}, and use the absorbing boundary condition at the left endpoint (representing underground region) and stress-free boundary condition at the right endpoint (representing earth's surface).
All the rest terms are randomly generated from the same distribution as explained in Section~\ref{sec:wave}, including choosing $Lu$ randomly
from $-c^2(x)u_{xx},-c(x)(c(x)u_x)_x,-(c^2(x)u_x)_x$ with equal probability.
For each equation, we assume access to solutions under 32 different initial conditions $g(x),h(x)$ and source terms $s_T(t)s_X(x)$, with observable samples limited to 32 fixed spatial locations and an average of 32 accessible samples per location.
We add noise with a relative amplitude of $1\%$ to the observed values at these points and $3\%$ to the initial values used as inputs for PDEformer-1.
The results are shown in Figure~\ref{fig:70}.

\begin{figure}[htpb]
	\centering
	\includegraphics[width=\linewidth]{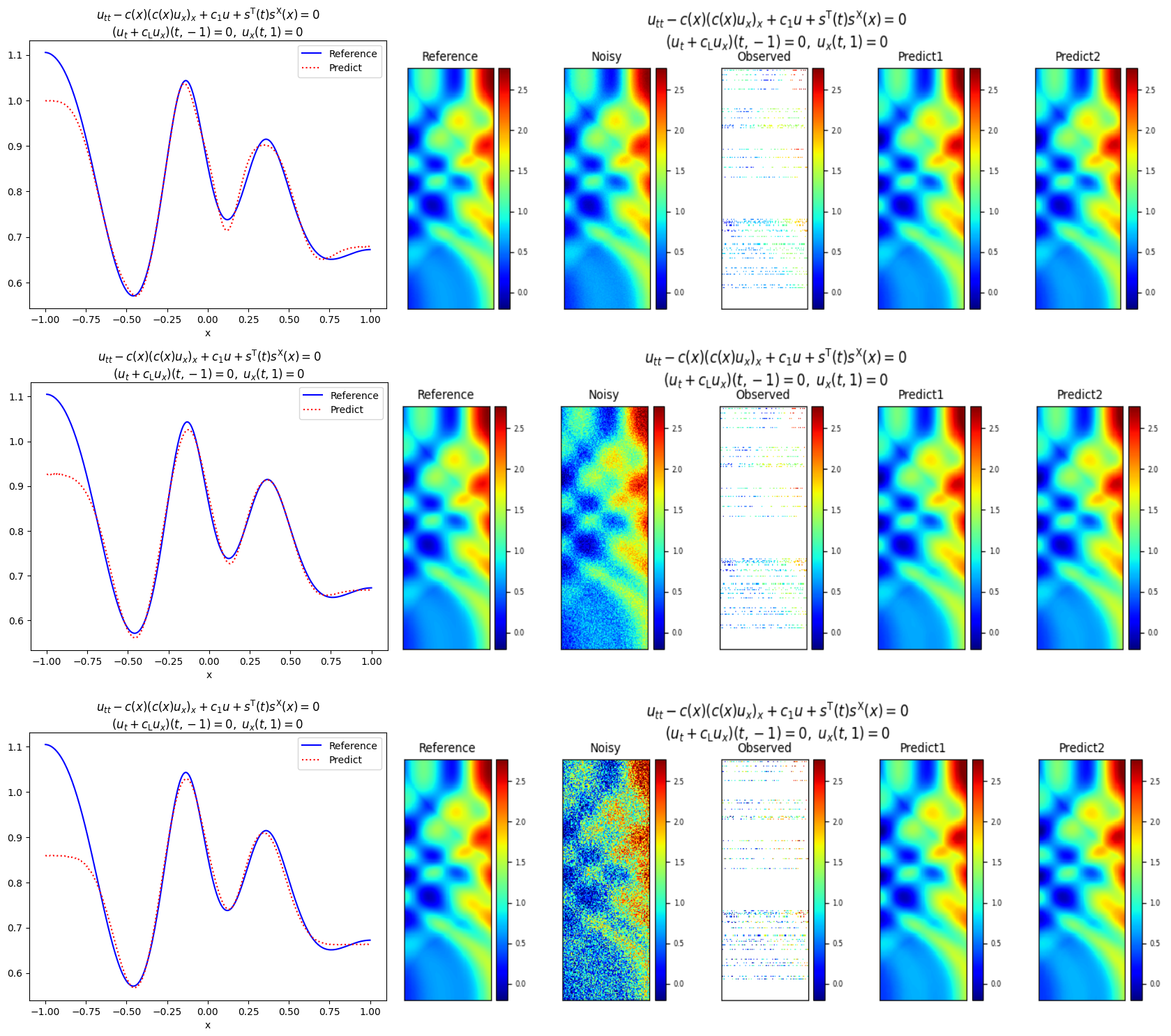}
	\caption{Observed solutions at different noise levels (taking one of the initial conditions as an example), along with the corresponding recovery results for the wave velocity field. The predicted velocity field is plotted with dashed lines, and the true value is plotted with solid lines. From the top to the bottom, we add noise with a relative amplitude of 1\%, 10\%, 30\% to observed values and 3\%, 10\%, 20\% to initial conditions, respectively.}%
	\label{fig:19}
\end{figure}

We have also conducted experiments with different noise levels, with the examples of observed solutions and the recovery results shown in Figure~\ref{fig:19}.

\section{Discussions}
\label{sec:discussion}
PDEformer-1 has demonstrated strong capabilities in several aspects.
In terms of versatility, unlike common neural operators, PDEformer-1 expresses the symbolic form of PDEs as a computational graph and uses it as input to the network, allowing to solve different forms of equations using a single model.
Even considering only the numeric information in PDEs, the inputs (initial values, coefficient fields, etc.) and outputs (predicted solutions of the equations) of PDEformer-1 are not dependent on a specific set of discrete grids.
Although the experiments in this paper focus on solving one-dimensional time-dependent PDEs with a single unknown variable, the PDEformer-1 architecture is capable of handling a much more comprehensive range of equation types.
We can use a single model to solve equations of different spatial dimensions, both time-dependent and time-independent equations, and systems of PDEs with varying numbers of variables and equations.
This makes it more versatile than any existing architecture known to the authors.

In terms of solution accuracy, PDEformer-1 also has a significant advantage over specialized models that can only be trained for specific equations.
For PDEs within the distribution of the pretraining data, the accuracy of PDEformer-1's direct inference exceeds that of specifically trained models.
If the coefficients of a PDE are beyond the range of the pretraining data, PDEformer-1's direct inference accuracy may not be ideal, but with a small amount of data for quick finetuning, its accuracy can still surpass that of specialized models.
Even when solving equations with previously unseen forms (such as wave equations), PDEformer-1 can perform incremental training based on the existing training results without the need of learning from-scratch.
Additionally, using inverse problems (recovering scalar coefficients of the equations, system identification, and recovering source fields and wave velocity fields) as examples, we primarily demonstrate the strong potential of a trained PDEformer-1 for various downstream tasks.

It is necessary to acknowledge that the current PDEformer-1 model has only been trained on one-dimensional time-dependent equations, mainly as a proof of concept, and its practical value is still relatively limited.
For practical problems in scientific research, engineering, and medical fields, the focus is primarily on two-dimensional, three-dimensional, or even higher-dimensional equations.
To apply PDEformer to these real-world scenarios in the future, we will explore efficient training methods for PDEformer so that it can solve two-dimensional and three-dimensional equations, as well as how we can facilitate data from equations of different dimensions to benefit our model.
The range of equation types examined will also be further broadened, potentially involving diverse algebraic or differential operations, more special functions, and complex domain shapes that may appear in two-dimensional and three-dimensional PDEs.
Additionally, we will explore more types of downstream tasks.
Beyond the inverse problems preliminarily explored in this paper, we will apply PDEformer to other tasks such as optimization design, control, and uncertainty quantification in the future.

\section*{Acknowledgments}
This work is supported in part by the National Science and Technology Major Project (2022ZD0117804).
Bin Dong is supported in part by the New Cornerstone Investigator Program.

\bibliographystyle{plainnat}

\bibliography{Ref}
\appendix
\section*{Appendix}
In the Appendix, we offer comprehensive supplementary materials to enhance the understanding of our study and support the reproducibility of our results. Appendix~\ref{sec:PDEDAGdetails} delves into the details of our computational graph representation, elucidating its design and the rationale behind its structure. 
Appendix~\ref{sec:arch} presents the details of PDEformer-1 and all the baseline models, along with the training settings.
Appendix~\ref{sec:PDEBenchData} explains the PDEBench datasets' features and our postprocessing steps.
\section{Details of the Computational Graph Representation}
\label{sec:PDEDAGdetails}
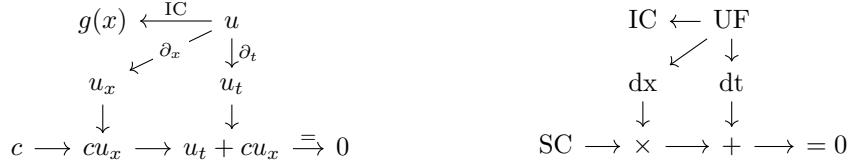
\begin{figure}[htbp]
	\centering
	\begin{minipage}[p]{0.5\textwidth}
		\begin{tikzcd}[row sep=small, column sep=small]
					& g(x)           & u \arrow[d, "\partial_t"] \arrow[l, "\mathrm{IC}"'] \arrow[ld, "\partial_x" description] &   \\
					& u_x \arrow[d]  & u_t \arrow[d]                                                                            &   \\
		c \arrow[r] & cu_x \arrow[r] & u_t+cu_x \arrow[r, "="]                                                                  & 0
		\end{tikzcd}
	\end{minipage}%
	\begin{minipage}[p]{0.4\textwidth}
		\begin{tikzcd}[row sep=small, column sep=small]
							  & \mathrm{IC}           & \mathrm{UF} \arrow[d] \arrow[ld] \arrow[l] &    \\
							  & \mathrm{dx} \arrow[d] & \mathrm{dt} \arrow[d]                      &    \\
		\mathrm{SC} \arrow[r] & \times \arrow[r]      & + \arrow[r]                                & =0
		\end{tikzcd}
	\end{minipage}%
	\caption{Illustration of how the form of a PDE can be represented as a computational graph, taking the advection equation $u_t+cu_x=0,u(0,x)=g(x)$ as the example. 
		The left panel shows the logical meaning of the nodes and edges, and the right panel illustrates the formalized data structure that is taken as the input of PDEformer-1.
		We also note that different from textual representations, this formalization of DAG is independent of the choice of symbols and the order of addition or multiplication.
		For example, the equation $\beta v_x+v_t=0,v|_{t=0}=v_0$ also corresponds to the DAG shown on the right panel.
	}
	\label{figs:PDEasDAG}
\end{figure}
Figure~\ref{figs:PDEasDAG} gives an illustration of the semantic meanings of the computational graph.
We also have the following remarks on the computational graph representation of the PDE form:
\begin{itemize}
	\item Only a small number of node types are involved in this computational graph:
		\verb|UF| (unknown field variable), \verb|SC| (scalar coefficient), \verb`CF` (coefficient field), \verb`VC` (varying coefficient), \verb|IC| (initial condition), $\mathtt{|x_L}$, $\mathtt{|x_R}$ (boundary value on the left/right endpoint),
		$\mathrm{dt,dx}$ (differentiation with respect to $t$ and $x$, respectively),
		$+$ (sum), $\times$ (product), $-$ (negation), $(\cdot)^2$ (square), $\sin,\cos$ (trigonometric functions), $=0$ (being equal to zero).
	\item Note that `$-$' stands for negation rather than subtraction in the computational graph, making it a unary rather than a binary operation.
		This is because nodes in this computational graph are not ordered, and the edges are homogeneous. 
		If the binary subtraction operation is involved in the computational graph, we cannot differentiate the subtrahend and the minuend.
	\item For the same reasons, we do not include a power operation node. However, we note that powers with a positive integer exponent can be expressed.
		For example, $u^3=u\times u^2$, and $u^{11}=((u^2)^2)^2\times u^2\times u$ since $11=2^3+2^1+2^0$.
	\item Although not involved in our experiments, node types representing other special functions such as $\exp,\log$ can be introduced as well, depending on the form of the PDEs involved.
		This also enables expression of the general power operation, since we have $a^b=\exp(b\times\log(a))$.
	\item 
        Disregarding the auxiliary nodes, nodes with type \verb|UF|, \verb|SC|, \verb|CF| and \verb|VC| would have a zero in-degree, $+$ and $\times$ have an in-degree that is equal to or greater than two, and the in-degrees of all the remaining nodes would be exactly one.
	\item In terms of the auxiliary nodes, we let each branch node $\mathtt{b}_k$ receive an edge from the corresponding \verb|IC| node, or emanate an edge towards the \verb`CF` or \verb`VC` node, and each latent modulation node $\mathtt{m}_\ell$ emanate an edge towards \verb|UF|.
		We adopt such a convention of edge direction in order to improve the connectivity of the final DAG, since we shall mask out the attention between disconnected node pairs in the graph Transformer module (see Appendix~\ref{sec:detail_DAG}).
\end{itemize}
\subsection{Detailed Generation of the Input Feature Vector} 
\label{sec:fun_enc}
We shall explain in further detail how the input feature embedding vector $\xi_i\in\R^{d_e}$ of node $i$ is generated.
If node $i$ is of type \verb|SC|, $\mathtt{|x_L}$ or $\mathtt{|x_R}$ with the corresponding scalar value $c\in\R$, we shall take
$\xi_i=\text{Scalar-Encoder}(c),$
where $\text{Scalar-Encoder}$ is a multi-layer perceptron (MLP) that has two hidden layers with $256$ neurons each.
Alternatively, if node $i$ is of type $\mathtt{b}_k$ with $k\in\{1,\dots,N\}$, and the corresponding function is $s(x)$, we take instead
$\xi_i=\text{Function-Encoder}(\{(x_j,s(x_j))\mid j\in \mathcal{J}\})_k.$
Inspired by VIDON~\cite{VIDON}, the function encoder is designed as an adapted variant of the DeepSet~\cite{DeepSet} network, in the form
\[\begin{split}
	&\qquad\text{Function-Encoder}(\{(x_j,s(x_j))\mid j\in \mathcal{J}\})
	\\&= \psi^3\left(\sum_{j\in \mathcal{J}}\mathrm{softmax}_j\Bigl(\psi^2(x_j,s(x_j))\Bigr)\odot\psi^1(x_j,s(x_j))\right)
	\\&= \psi^3\left(\frac{\sum_{j\in \mathcal{J}}\exp\Bigl(\psi^2(x_j,s(x_j))\Bigr)\odot\psi^1(x_j,s(x_j))}{\sum_{j\in \mathcal{J}}\exp\Bigl(\psi^2(x_j,s(x_j))\Bigr)}\right)
,\end{split}\]
where $\psi^1,\psi^2:\R^2\to\R^{512}$ and $\psi^3:\R^{512}\to\R^{N\times d_e}$ are all represented by MLPs with two hidden layers of width $d_\psi=512$, the exponential and division operations are applied elementwise, and $\odot$ denotes elementwise multiplication of two vectors.
We observe that taking $N=4$ and $\#\mathcal{J}=256$ already produces satisfactory accuracy, and stick to this setting throughout the experiments.
For all the remaining cases, the input feature vector is taken as $\xi_i=\text{Scalar-Encoder}(0)$.
\subsection{Computational Graph of Wave Equation}
\begin{figure}[ht]
  \centering
  \includegraphics[width=0.6\linewidth]{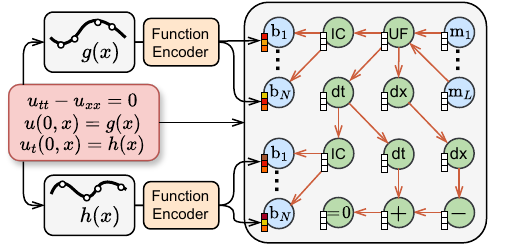}
  \caption{Computational graph of the wave equation, containing $\texttt{dt} \to \texttt{dt}$.}
  \label{fig:8}
  \end{figure}

We illustrate the construction of the computational graph of a simple wave equation in Figure~\ref{fig:8}.
Different from equations in PDEBench, wave equations include elements not seen by PDEformer-1 during the pretraining stage (one $dt$ connected to another $dt$ in the graph).
Although the actual wave equation involved in our experiments has a more complicated form~\eqref{eq:wave}, the way of representing the term $u_{tt}$ is the same as in this simple example.

\subsection{Boundary Conditions}
\begin{figure}[ht]
\centering
    \includegraphics[width=0.5\linewidth]{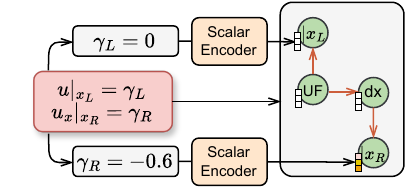}
    \caption{example of computational graph of PDE with non-periodic boundary condition, with Dirichlet boundary condition at the left endpoint and Neumann boundary condition at the right endpoint.}
    \label{fig:DirNeum}
\end{figure}

\begin{figure}[ht]
    \centering
    \includegraphics[width=0.5\linewidth]{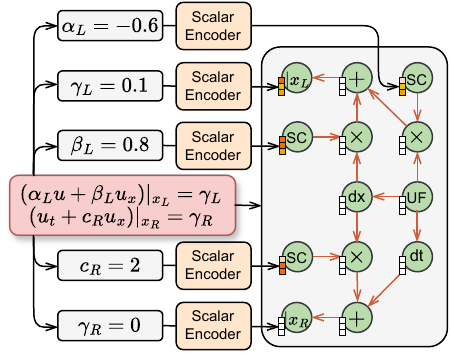}
    \caption{example of computational graph of PDE with non-periodic boundary condition, with Robin boundary condition at the left endpoint and Mur boundary condition at the right endpoint.}
    \label{fig:RobinMur}
\end{figure}

If the PDE uses periodic boundary conditions, there is no need to introduce additional nodes in the corresponding computational graph.
For non-periodic boundaries, we shall introduce two node types $ \mathtt{|x_L},\mathtt{|x_R}$ to specify the boundary values at the left and the right endpoints of the spatial domain.
Figure~\ref{fig:DirNeum} and~\ref{fig:RobinMur} show two examples of how different boundary conditions can be specified in the computational graph.
For example, the connection of the \verb|UF| node to $\mathtt{|x_L}$ in Figure~\ref{fig:DirNeum} represents $u|x_L$, i.e. the value of $u$ at the boundary point $x=x_L$.
Similarly, the $\texttt{UF}\to\texttt{dx}\to\mathtt{|x_R}$ element represents the value of $u_x$ at $x=x_R$.
Additional algebraic and differential operations have to be introduced to specify more complicated boundary conditions, as can be seen in Figure~\ref{fig:RobinMur}.
\subsection{Non-Constant Coefficients}
\begin{figure}[ht]
    \centering
    \includegraphics[width=0.5\linewidth]{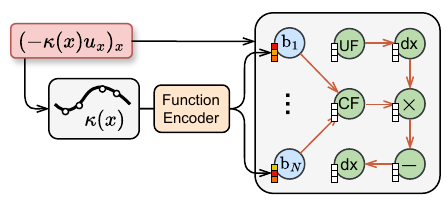}
    \caption{example of a computational graph including a coefficient field.}
    \label{fig:PDEformer-CF}
\end{figure}

A coefficient field $s(x)$ with spatial dependency is represented by a node with type \verb|CF|.
Figure~\ref{fig:PDEformer-CF} illustrates the computational graph of a diffusion term with a non-constant diffusivity coefficient $\kappa(x)$.

A varying coefficient with temporal dependency (i.e., of the form $s(t)$, which is a function of time $t$) can be represented in a similar way by a \verb|VC| node.
Moreover, a varying coefficient field of the form $s^\mathrm{T}(t)s^\mathrm{X}(x)$ (assumed to be separable in time and space) can be expressed in the computational graph with the help of an additional multiplication operation.
\subsection{Multi-Variable Equations}
\begin{figure}[htbp]
    \centering
    \includegraphics[width=0.6\linewidth]{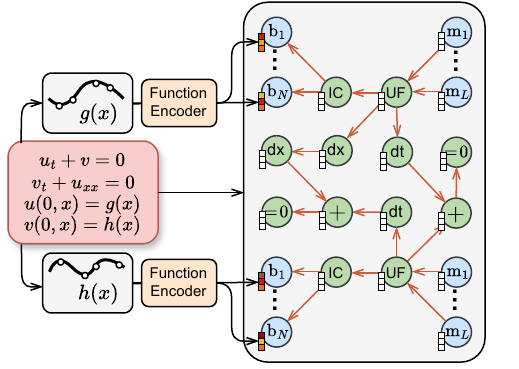}
    \caption{An example of the computational graph of a multi-variable PDE. Two different nodes of type \texttt{UF} have been introduced to represent the two variables $u$ and $v$, respectively.}
    \label{fig:mCompn}
\end{figure}

Although all the experiments in this study involve only PDEs with a single variable $u$, we note that
systems of partial differential equations with multiple variables (or, in other words, multiple components of the equation solution) may also be represented using computational graphs and can therefore be solved by PDEformer-1 in principle.
A simple example is given in Figure~\ref{fig:mCompn}, where we have introduced one node of type \verb|UF| for each unknown PDE variable involved.
\subsection{Efficient Implementation of Input Feature Vector Assignment}
\begin{figure}[htbp]
    \centering
    \includegraphics[width=1.0\linewidth]{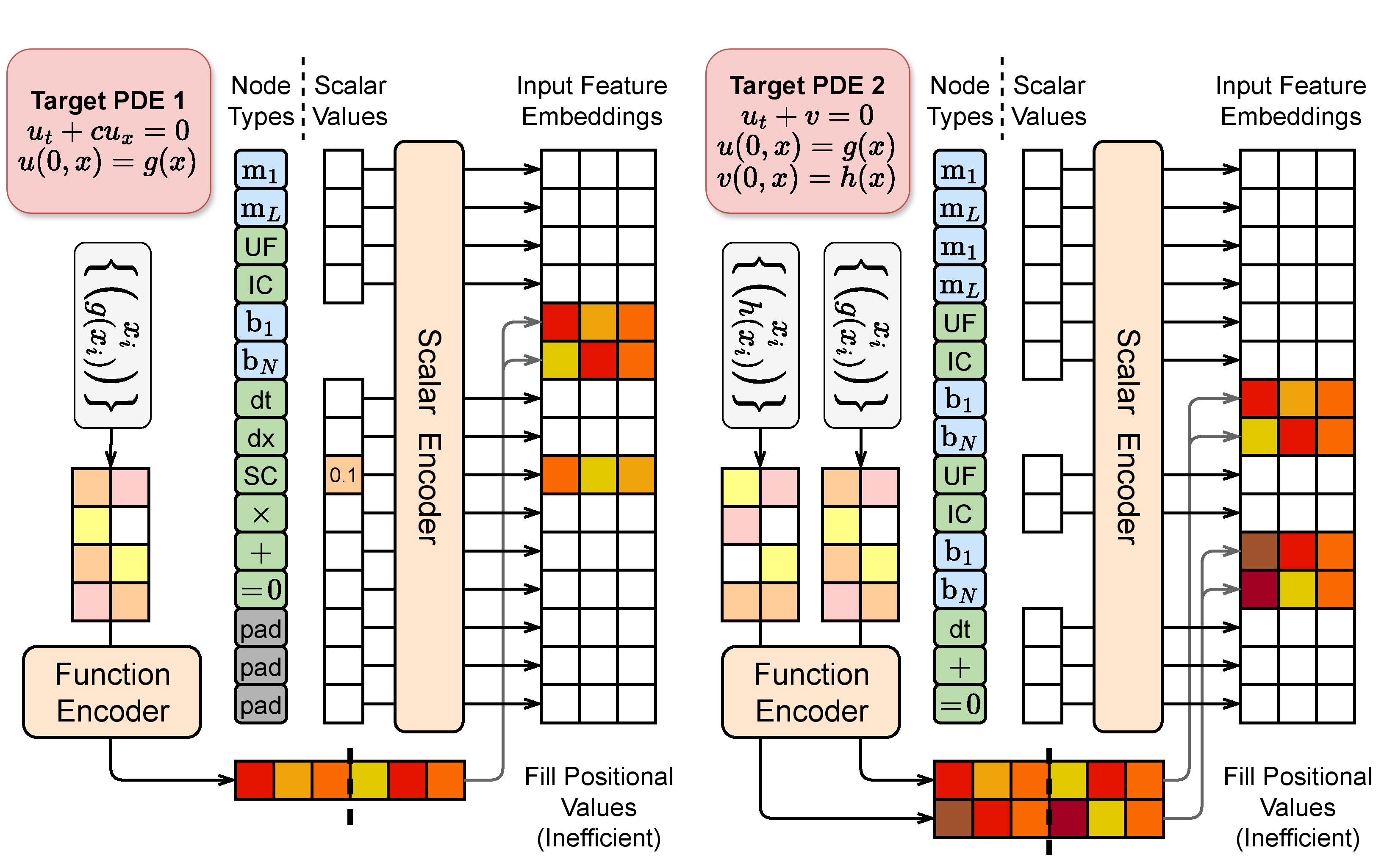}
    \caption{The na\"ive method of input feature vector assignment.}
    \label{fig:FuncEmb-naive}
\end{figure}

In the left part of Figure~\ref{fig:PDEformerV2Arch}, we assign an input feature embedding of dimension $d_e$ to each node in the graph.
Some nodes receive input features from scalar encoders, while others receive them from function encoders.
For this setup, the most direct and na\"ive method of implementation is shown in Figure~\ref{fig:FuncEmb-naive} (for simplicity of illustration, we set $L=N=2$ and use an incomplete equation for PDE 2).

For this approach, we need to examine the information of each input function involved in the PDE, determine the permutation position of the corresponding `branch' nodes ($\mathtt{b}_1,\dots,\mathtt{b}_N$) among all graph nodes, and fill the output of the function encoder into the corresponding position of the overall input feature tensor%
\footnote{In programming implementation, all input and output data of the network are represented in tensor format.}%
containing all graph nodes.
In the training process, a data batch often contains multiple different PDEs, each with a different number of input functions, and the permutation positions of the corresponding branch nodes are also different.
If this implementation method is used, the code will not only be cumbersome to write but also inefficient at runtime.
\begin{figure}
    \centering
    \includegraphics[width=1.0\linewidth]{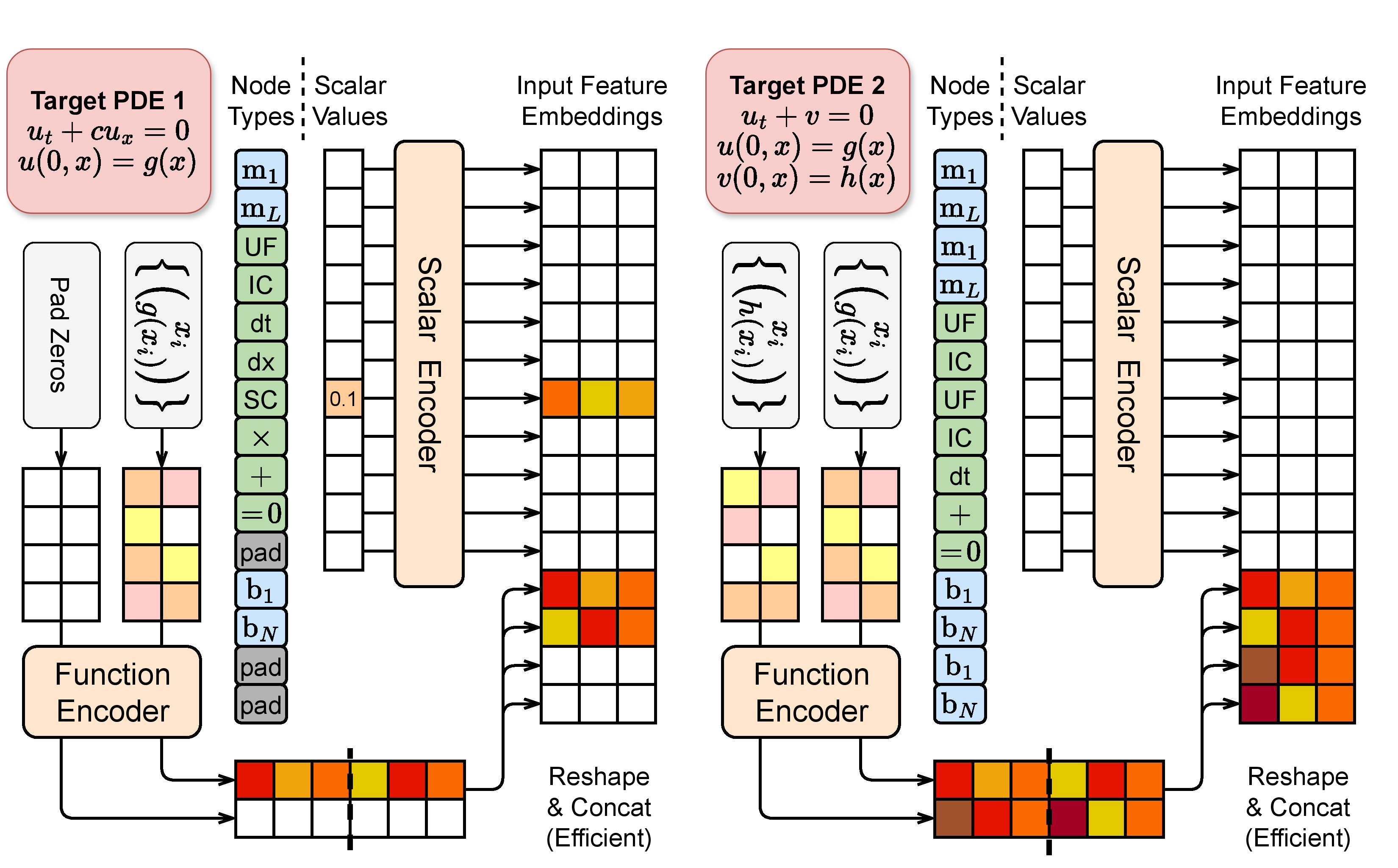}
    \caption{How input feature vector assignment is implemented in practice.}
    \label{fig:FuncEmb-actual}
\end{figure}

In programming implementation, we adopt the following alternative approach, which is shown in Figure~\ref{fig:FuncEmb-actual}.
\footnote{Here, we take advantage of the special attribute possessed by the graph Transformer.
Unlike most Transformers designed for sequence data, the forward computation process of the graph Transformer only makes use of the topological information of the graph, and does not rely on the order of nodes.
This allows us to rearrange the nodes in an arbitrary order.}
This approach specifies that the input features of nodes with earlier positions are generated by scalar encoders, while those with later positions are given by function encoders.
With the relatively simple tensor shape modifications (Reshape) and tensor concatenation (Concat) operations, the final overall input feature tensor can be obtained without the need for tedious judgment and value assignment processes.

\section{Model Architecture of Our Model and Baseline Models}
\label{sec:arch}
\subsection{Details of the Graphormer Structure}
\label{sec:detail_DAG}
\paragraph{Initial Embedding Vector} In the graph Transformer, the initial embedding of node $i$ is given as:
\[
    h_i^{(0)} = x_{\text{type}(i)} + \xi_i + z^-_{\text{deg}^-(i)} + z^+_{\text{deg}^+(i)}
\]
where $\xi_i\in\R^{d_e}$ is the input feature embedding vector of node $i$, and $x,z^-,z^+\in\R^{d_e}$ are learnable vectors specified by node type $\text{type}(i)$,  in-degree $\text{deg}^-(i)$ and out-degree $\text{deg}^+(i)$, respectively.

\paragraph{Attention Bias} We denote $\phi(i,j)$ to be the shortest path length from node $i$ to node $j$. If such a path does not exist or has a length greater than $14$, we set $\phi(i,j) = 14$. For each attention head involved in the graph Transformer, the attention bias according to the node pair $(i,j)$ is given as 
\[
    B_{ij} = b^+_{\phi(i,j)} + b^-_{\phi(j,i)} + d_ij
\]
Here, $b^+_{\phi(i,j)},b^-_{\phi(i,j)}$ are learnable scalars indexed by $\phi(i,j)$ and $\phi(j,i)$ respectively, and is shared across all layers. The additional term $d_{ij}$, which does not appear in the original Graphormer, is introduced to mask out attention between disconnected node pairs.
More specifically, when node $i$ and node $j$ are connected in the graph, i.e. there exists a path either from $i$ to $j$ or from $j$ to $i$, we take $d_{ij}=0$, and set $d_{ij}=-\infty$ otherwise. We observe in our experiments that the overall prediction accuracy can be improved with such an additional masking operation. Moreover, since our graph has homogeneous edges, we do not introduce the edge encoding term that appears in the original Graphormer.

\paragraph{Graph Transformer Layer} The structure of the graph transformer layer is the same as the original Graphormer, and we include it here for convenience to the readers.
Each layer takes the form
\begin{align}
    \bar h^{(l)} &= \text{Attn}(\text{LN}(h^{(l-1)})) + h^{(l-1)}\\
    h^{(l)} &= \text{FFN}(\text{LN}(\bar h^{(l)})) + \bar h^{(l)}
,\end{align}
where $\text{FFN}$ represents a position-wise feed-forward network with a single hidden layer and GeLU activation function, and $\text{LN}$ stands for layer normalization.
In terms of the self-attention block $\text{Attn}$, we shall follow the convention in the original Graphormer paper and only present the single-head case for simplicity.
Let $H = [h_1', \cdots, h_n']^\mathrm{T}\in\R^{n\times d_e}$ denote the input of the self-attention module involving $n$ graph nodes, the self-attention is computed as
\[
\begin{aligned}
    Q = HW_Q,\quad K = HW_K,\quad V = HW_V,\\
    A = \frac{QK^\mathrm{T}}{\sqrt{d_e}} + B, \quad \text{Attn}(H) = \text{softmax}(A)V,
\end{aligned}
\]
where $W_Q,W_K,W_V\in\R^{d_e\times d_e}$ are the projection matrices, and $B$ is the attention bias given before.
The extension to the multi-head attention is standard and straightforward.  
\subsection{Details of the INR Structure}
\label{sec:inr}
In the realm of Implicit Neural Representation (INR), data samples are interpreted as coordinate-based functions, where each function accepts a coordinate $(t,x)$ as input and yields an approximated function value $\hat{u}(t,x)$ at that specific coordinate point.
Various architectures of such INRs have been proposed in the literature, including DeepONet~\citep{DeepONet}, HyperDeepONet~\citep{HyperDeepONet} for neural operators, as well as SIREN~\citep{SIREN}, WIRE~\citep{WIRE}, MFN~\citep{MFN}, Poly-INR~\citep{PolyINR} and others~\citep{GaussINR,TransINR,ShapE} in computer vision.
In the experiments, we utilize an adapted version of Poly-INR~\citep{PolyINR}, which exhibits better prediction accuracy and training stability compared with other candidates in our setting.
Inspired by COIN++~\citep{coinpp}, we also make a modification to the $L$ hypernets in the Poly-INR structure, in which the $\ell$-th hypernet takes $\mu^\ell\in\R^{d_e}$ as its input, and generates the scale- and shift-modulations for the $\ell$-th hidden layer of our Poly-INR.
\begin{figure}[ht]
  \centering
  \includegraphics[width=\textwidth]{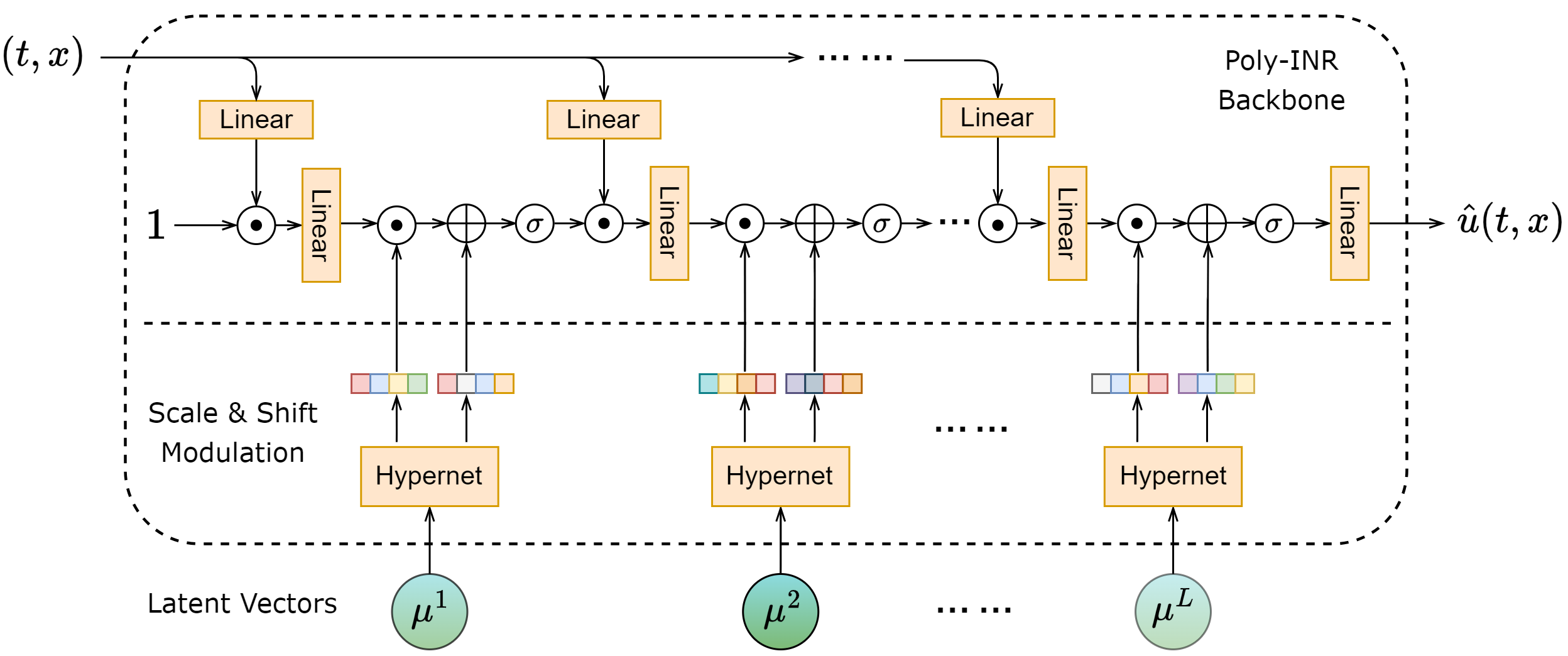}
  \caption{INR decoder architecture of PDEformer.}
  \label{fig:INR-decoder}
  \end{figure}

The intricate architecture of our INR decoder is illustrated in figure~\ref{fig:INR-decoder}, with the mathematical framework detailed below.
We take $h_0=\mathbf{1}$ to be the vector with all entries equal to one. For $\ell=1,2,\dots,L$, we compute
\[
\begin{aligned}
	g_\ell &= W^{\text{in}}_{\ell} \begin{bmatrix} t \\ x \end{bmatrix} + b^\text{in}_\ell, \quad
	s_{\ell}^{\text{scale}} = \text{MLP}_\ell^{\text{scale}}(\mu^\ell), \quad
	s_{\ell}^{\text{shift}} = \text{MLP}_\ell^{\text{shift}}(\mu^\ell), \\
	q_{\ell} &= s_{\ell}^{\text{scale}} \odot \left( W^{\text{h}}_{\ell} \left(h_{\ell-1} \odot g_\ell\right) + b_\ell^\text{h} \right) + s_{\ell}^{\text{shift}}, \quad
	h_\ell = \sigma\left(q_{\ell}\right),
\end{aligned}
\]
and the network output is given as $\hat{u}(t,x) = W^{\text{Last}}h_{L} + b^{\text{Last}}$.
Here, the activation function $\sigma(\cdot)$ is a leaky-ReLU operation with a slope of $0.2$ at the negative input range, followed by a clipping operation into the interval $[-256,256]$ to improve training stability.
The hypernets correspond to $\text{MLP}_\ell^{\text{scale}}$ and $\text{MLP}_\ell^{\text{shift}}$.
Note that in the original Poly-INR, the hypernets are utilized to generate $W^\text{in}_\ell$ and $b^\text{in}_\ell$.
Compared with our practice of generating $s^\text{scale}_\ell$ and $s^\text{shift}_\ell$, this method exhibits better accuracy, but deteriorates the training efficiency, and is therefore not adopted in our experiments.
\subsection{Training Setting}
\label{sec:setting}

The experimental settings, including model hyperparameters and configurations, are outlined in Table \ref{tab:train_settings}. 
For a comprehensive understanding of the baseline models employed in our experiments, we provide an overview of all models:

\begin{itemize}
    \item \textbf{DeepONet:} DeepONet employs a unique architecture with two sub-networks: a branch net and a trunk net. The branch net processes a fixed number of sensor observations ($256$ points from the initial condition in our case), while the trunk net handles coordinate inputs for inference, akin to PDEformer-1's input mechanism. The outputs from both networks are combined to produce the solution value. Each sub-network consists of a six-layer MLP with $256$ hidden neurons and utilizes the ReLU activation function. Notably, DeepONet's mesh-free nature allows for training with scattered data points, enabling us to sample $8192$ points per iteration from $256\times256$ grids for each data sample during both DeepONet's training and PDEformer-1's finetuning processes.
    
    \item \textbf{FNO:} The Fourier Neural Operator (FNO) operates on a mesh-dependent yet resolution-independent principle. It initially transforms regular grid data into multi-channel hidden features through a pointwise fully connected layer, followed by processing through several Fourier Layers, and finally map to the solution grid. In Fourier Layer, the FNO keeps the lowest $12$ Fourier modes. In our experiments, the FNO2D model is utilized, with the initial condition ($256$ spatial points) extended to form a $256\times256$ input grid, allowing for simultaneous full field output.
    
    \item \textbf{U-Net:} U-Net adopts a CNN-based encoder-decoder framework, distinguished by its 4 layers of downsampling and upsampling convolutions, bridged by intermediate residual connections. Analogous to FNO2D, both the input and output dimensions are set to $256\times256$. Unlike the mesh-free DeepONet or PDEformer-1, FNO and U-Net require training data organized in regular grids.
    
    \item \textbf{PDEformer-1:} The Transformer-based Graphormer is configured with $9$ layers, a $512$-dimensional embedding space, and $32$ attention heads. The Poly-INR part employs $L=8$ hidden layers with $256$ neurons, and each hidden layer is dynamically modulated using separate scale and shift hypernets, each comprising of a $3$-layer MLP with independent parameters.
\end{itemize}

In the pretraining stage of PDEformer-1, we employ nRMSE loss function due to its effectiveness in improving training efficiency. The nRMSE loss function can be mathematically represented as:
\[
    \mathcal{L}_{\text{nRMSE}} = \frac{\|u-\hat{u}\|_{L_2}}{\|u\|_{L_2}}
\]
A learning rate schedule is implemented, progressively reducing the learning rate at predetermined epochs to improve the stability of the training process. Moreover, a warm-up period is utilized at the start of training to mitigate the risk of early training failures by gradually increasing the learning rate from zero to the initial pre-scheduled value.
\begin{table}[ht]
    \centering
    \caption{Hyperparameters. We give out the differences of PDEformer-1 with scale S, M, L, and XL. By $a,b,c,d$ in the `value' part, we mean that the value of the corresponding hyperparameter is $a$ for the S-scale model, $b$ for the M-scale model, and analogously for $c$ (L) and $d$ (XL).}
    \begin{tabular}{lll}
    \hline
    \multicolumn{1}{c}{\textbf{Parameter}} & \multicolumn{1}{c}{\textbf{Value}} & \multicolumn{1}{c}{\textbf{Description}} \\
    \hline
    \multicolumn{3}{l}{\textbf{DeepONet}} 
    \\
    \quad trunk\_dim\_in & 2 & Input dimension of the trunk network \\
    \quad trunk\_dim\_hidden & 256 & Dimension of hidden features in the trunk network \\
    \quad trunk\_num\_layers & 6 & Number of layers in the trunk network \\
    \quad branch\_dim\_in & 256 & Input dimension of the branch network\\
    \quad branch\_dim\_hidden & 256 & Dimension of hidden features \\
    \quad branch\_num\_layers & 6 & Number of layers in the branch network \\
    \quad dim\_out & 2048 & Output dimension of the trunk net and the branch net \\
    \quad num\_tx\_samp\_pts & 8192 & Number of sample points used per training iteration \\
    \quad learning\_rate & \(0.0003\) & The initial learning rate for the optimizer \\
    \multicolumn{3}{l}{\textbf{FNO}} \\
    \quad resolution & 256 & The resolution of the grid \\
    \quad modes & 12 & The truncation number of Fourier modes\\
    \quad channels & 20 & The number of channels in the hidden layers \\
    \quad depths & 4 & The number of Fourier Layers in the neural network \\
    \quad learning\_rate & \(0.0001\) & The initial learning rate for the optimizer \\
    \multicolumn{3}{l}{\textbf{U-Net}} \\
    \quad learning\_rate & \(0.0001\) & The initial learning rate for the optimizer \\
    \multicolumn{3}{l}{\textbf{Autoregressive U-Net}} \\
    \quad learning\_rate & \(0.0001\) & The initial learning rate for the optimizer \\
    \multicolumn{3}{l}{\textbf{PDEformer-1}} 
    \\
    \multicolumn{3}{l}{\textbf{\quad Scalar Encoder}} 
    \\
    \qquad dim\_hidden & 256 & Dimension of the hidden feature\\
    \qquad num\_layers & 3 & Number of hidden layers\\
    \multicolumn{3}{l}{\textbf{\quad Function Encoder}}
    \\
    \qquad num\_branch & 4 & Number of branches nodes $N$ used for each coefficient field $s(x)$ \\
    \qquad dim\_hidden & 512 & Dimension of the hidden feature\\
    \qquad num\_layers & 6 & Number of hidden layers\\
    \multicolumn{3}{l}{\textbf{\quad Graphormer}} 
    \\
    \qquad num\_layers & 4, 6, 9, 12 & Number of layers in Graphormer \\
    \qquad embed\_dim & 128, 256, 512, 768 & Dimension of the feature embedding \\
    \qquad ffn\_embed\_dim & 128, 256, 512, 768 & Dimension of the feed-forward network embedding \\
    \qquad num\_heads & 16, 32, 32, 32 & Number of attention heads \\
    \qquad pre\_layernorm & True & Whether to use layer normalization before each block \\
    \multicolumn{3}{l}{\textbf{\quad Poly-INR}} \\
    \qquad dim\_in & 2 & Input dimension \\
    \qquad dim\_hidden & 64, 128, 256, 512 & Dimension of the hidden feature \\
    \qquad dim\_out & 1 & Output dimension \\
    \qquad num\_layers & 4, 6, 9, 12 & Number of hidden layers \\
    \multicolumn{3}{l}{\textbf{\quad Layerwise Hypernet}} \\
    \qquad hyper\_dim\_hidden & 256 & Dimension of hidden layers in a hypernet \\
    \qquad hyper\_num\_layers & 3 & Number of layers in a hypernet \\
    \qquad share\_hyper & False & Whether hypernets share parameters across all layers \\ 
    \multicolumn{3}{l}{\textbf{PDEformer-1 Pretraining}} \\
    \quad batch\_size & 640 & Total batchsize used in one iteration \\
    \quad learning\_rate & \(0.0001\) & The initial learning rate for the optimizer \\
    \quad epochs & 500 & The total number of training epochs \\
    \quad loss\_type & nRMSE & Use the normalized root-mean-squared-error for training \\
    \quad optimizer & Adam & The optimization algorithm \\
    \quad lr\_scheduler & cosine & The learning rate scheduler \\
    \quad warmup\_epochs & 10 & Epochs to linearly increase the learning rate \\
    \hline
    \end{tabular}
    \label{tab:train_settings}
\end{table}
\section{PDEBench Datasets}
\label{sec:PDEBenchData}

In this subsection, we present an overview of three 1D PDE datasets derived from PDEBench, which we employed in our experimental analysis. Each dataset, tailored to a specific PDE type and coefficient configuration, encompasses 10k instances. For our training purposes, we utilized 9k samples from each dataset, reserving the remaining 1k samples for testing.

\begin{itemize}
	\item Burgers' equation%
		\footnote{The convection term is $\partial_x(u^2)$ rather than $\partial_x(u^2/2)$ due to an implementation issue in the PDEBench data generation code. See \url{https://github.com/pdebench/PDEBench/issues/51} for more details.}:
		$\partial_t u+\partial_x(u^2)=\frac{\nu}{\pi} \partial_{x x} u$ for $(t,x)\in[0,2]\times[-1,1]$, where $\nu \in \{0.1,0.01,0.001\}$.
		This equation is a fundamental partial differential equation from fluid mechanics.
	\item Advection equation: $\partial_t u+\beta \partial_x u=0$ for $(t,x)\in[0,2]\times[0,1]$, where $\beta \in \{0.1, 1\}$.
		The equation models the transport of a quantity $u$ without alteration in its form.
	\item Reaction-Diffusion equation: $\partial_t u=\nu \partial_{x x} u+\rho u(1-u)$ for $(t,x)\in[0,1]\times[0,1]$, where we only consider $\nu=1, \rho=1$.
		This equation represents a process combining chemical reaction and diffusion dynamics.
\end{itemize}

The initial conditions for each dataset are given by $u_0(x)=\sum_{k_i=k_1, \ldots, k_N} A_i \sin \left(k_i x+\phi_i\right)$, with frequency numbers $k_i = \frac{2\pi n_i}{L_x}$, where $n_i$ are integers randomly selected within a pre-determined range and $L_x$ is the length of the spatial domain, amplitudes $A_i$ are random numbers within $[0,1]$, and phases $\phi_i$ are chosen randomly from the interval $\left(0, 2\pi\right)$.
The absolute value function with a random signature, as well as restriction to a random sub-interval by multiplying a window function, are applied afterwards with $10\%$ probability each.
For the Reaction-Diffusion equation, the range of the initial condition is rescaled to the unit interval $[0,1]$.

In order to utilize the pretrained PDEformer-1 model to make predictions, we rescale the spatial-temporal coordinates to the range $(t',x')\in[0,1]\times[-1,1]$, and the resulting PDEs taken as the input of PDEformer-1 have the following form:
\begin{itemize}
	\item Burgers' equation: $\partial_{t'}u+\partial_{x'}(2u^2)-\frac{2\nu}{\pi}\partial_{x'x'}u=0$, where $t'=t/2,x'=x$.
	\item Advection equation: $\partial_{t'}u+\partial_{x'}(4\beta u)=0$, where $t'=t/2,x'=2x-1$.
	\item Reaction-Diffusion equation: $\partial_{t'} u-4\nu\partial_{x'x'}u+(-\rho)u+\rho u^2=0$, where $t'=t,x'=2x-1$.
\end{itemize}

To ensure equitable comparisons among the baseline models, we standardize the resolution of all PDEBench samples to $256\times256$.
More specifically, the original PDEBench datasets have a spatial resolution of $1024$, which is downsampled to $256$.
The original number of recorded time-steps is $201$ for the Burgers' and Advection datasets and $101$ for the one-dimensional Reaction-Diffusion dataset, and a linear interpolation is utilized to obtain a temporal resolution of $256$.
It is important to note that PDEformer-1 makes mesh-free predictions, enabling us to set the temporal resolution to $101$ for the pretraining dataset and $256$ for the PDEBench dataset.
For FNO and U-Net (non-autoregressive case), the initial value is repeated $256$ times to form the two-dimensional data with resolution $256\times 256$, and then taken as the network input.

\end{document}